\documentclass[]{interact}
\usepackage{epstopdf}% To incorporate .eps illustrations using PDFLaTeX, etc.
\usepackage{subfigure}% Support for small, `sub' figures and tables
\makeatletter
\def\BState{\State\hskip-\ALG@thistlm}
\makeatother
\usepackage{pgf, tikz}
\usetikzlibrary{calc, fit, arrows}
\tikzstyle{line}=[draw]
\tikzstyle{arrow}=[draw, -latex]
\usetikzlibrary{positioning}
\usepackage[linesnumbered, ruled, vlined, noend]{algorithm2e}
\usepackage{lineno}
\usepackage[colorlinks=true, allcolors=blue]{hyperref}
\usepackage{suffix}
\usepackage{amsmath,amssymb,amsfonts}
\usepackage{mathtools}% http://ctan.org/pkg/mathtools
\usepackage{graphicx}
\usepackage{textcomp}
\usepackage{xcolor}
\usepackage{float}
\usepackage{nccmath}
\usepackage{nicefrac}
\usepackage{listings}% http://ctan.org/pkg/listings
\usepackage{multirow}
\usepackage{booktabs}
\usepackage{siunitx}
\usepackage{mathrsfs}
\usepackage{array}
\usepackage{longtable}
\usepackage{amsthm}
\usepackage{dsfont}
\usepackage{xspace}
\usepackage{makecell}
\usepackage{soul}
\usepackage{geometry}
\usepackage{rotating}
\usepackage{multicol}

\usepackage{natbib}% Citation support using natbib.sty
\bibpunct[, ]{(}{)}{;}{a}{}{,}% Citation support using natbib.sty
\renewcommand\bibfont{\fontsize{10}{12}\selectfont}% Bibliography support using natbib.sty

\begin{document}

% \articletype{ARTICLE TEMPLATE}

\title{A Comparison of Different Approaches to Dynamic Origin-Destination Matrix Estimation in Urban Traffic}

\author{
\name{
    Nicklas Sindlev Andersen\textsuperscript{a}, Marco Chiarandini\textsuperscript{a} and Kristian Debrabant\textsuperscript{a}
    \thanks{Contact emails:     \{\texttt{{sindlev$|$marco$|$debrabant}}\}\texttt{@imada.sdu.dk}}
}
\affil{
    \textsuperscript{a}University of Southern Denmark (SDU), 55 Campusvej, Odense M, Denmark.
}
}

\maketitle

\begin{abstract}
    Given the counters of vehicles that traverse the roads of a traffic network,
    we reconstruct the travel demand that generated them expressed in terms of the
    number of origin-destination trips made by users. We model the problem as a
    bi-level optimization problem. 
    At the inner-level, given a tentative demand, we solve a Dynamic Traffic
    Assignment (DTA) problem to decide the routing of the users between their
    origins and destinations. Finally, we adjust the number of trips and their
    origins and destinations at the outer-level to minimize the discrepancy
    between the counters generated at the inner-level and the given vehicle counts
    measured by sensors in the traffic network.
    We solve the DTA problem by employing a mesoscopic model implemented by the
    traffic simulator SUMO. Thus, the outer problem becomes an optimization
    problem that minimizes a black-box Objective Function (OF) determined by the
    results of the simulation, which is a costly computation.
    We study different approaches to the outer-level problem categorized as
    gradient-based and derivative-free approaches. Among the gradient-based
    approaches, we look at an assignment matrix-based approach and an assignment
    matrix-free approach that uses the Simultaneous Perturbation Stochastic
    Approximation (SPSA) algorithm. Among the derivative-free approaches, we
    investigate Machine Learning (ML) algorithms to learn a model of the simulator
    that can then be used as a surrogate OF in the optimization problem.
    We compare these approaches computationally on an artificial network. The
    gradient-based approaches perform the best in terms of solution quality and
    computational requirements. In contrast, the results obtained by the ML
    approach are currently less satisfactory but provide an interesting avenue for
    future research.
\end{abstract}

\begin{keywords}
  Bi-level optimization;
  Dynamic origin-destination matrix estimation;
  Dynamic traffic assignment;
  Simulation-based optimization
\end{keywords}

\section{Introduction}\label{sec:introduction}

Traffic on the road networks of urban areas is continuously increasing.
Municipalities can act on the infrastructure to keep high mobility standards,
modifying it to avoid congestion and ensure reliability. Therefore, there is an
urge to provide decision-makers with tools to make informed decisions. For
this purpose, traffic simulation software has been developed to assess different
planning scenarios. To use these simulators, the \emph{network supply}, and the
\emph{travel demand} must be known with some certainty. The network supply is,
in general, defined as the maximum number of vehicles the network road
infrastructure can handle at a given time. In contrast, the demand can
generally be said to be the number of vehicles that would like to travel on the
network at a given time.

The network supply and travel demand interact dynamically, i.e., the users
travel in their vehicles from one location to another on the network and
interact with each other and with the road infrastructure. However, congestion
may occur, and users might change their decisions and re-distribute on the
network. If the demand and the network supply are known, it can be possible to
obtain accurate distributions of users on the given network by calculating an
equilibrium situation corresponding to solving a \emph{Dynamic Traffic
Assignment (DTA)} problem. Then, by varying the network supply, one can analyze
different scenarios to address a network design task. Consequently, if the
demand in a particular traffic network is unknown, it is desirable to estimate
it.

Travel demand can be modeled in many different ways but is commonly modeled
using \emph{dynamic Origin-Destination (OD) matrices}. These matrices describe
the aggregate demand pattern changing over time (dynamic) in a set of time
intervals of a much larger analysis period. Due to their simplicity and
conciseness, these dynamic OD matrices are usually the required input of many
traffic simulation software packages. Because of this, the problem of trying to
estimate dynamic OD matrices has gained a lot of attention from researchers in
recent years.

The problem of estimating dynamic OD matrices is usually referred to as the
\emph{Dynamic OD Estimation (DODE)} problem in the literature. Different
formulations exist, but it is commonly formulated as a bi-level optimization
problem consisting of an inner and an outer optimization problem. A proposed
dynamic OD matrix is given as input to the inner problem, and a DTA problem is
solved. The result is a routing of the users in the network between origins and
destinations at every time interval. In the outer optimization problem, the
dynamic OD matrices are adjusted to minimize the discrepancy between the routing
returned by the solution of the inner problem and the one observed in real-life,
i.e., the number of vehicles that transited in every time interval on every road
of the network as measured by sensors. The observed routing is usually expressed
in terms of the number of vehicles, average speeds, densities, etc., that have
been captured by sensors located along roads in the traffic network under
investigation. However, other kinds of traffic observations can also be included
(e.g., travel times, intersection turning ratios, etc.). This information is
assumed to be directly related to the unknown travel demand and can guide the
retrieval of such a demand. In other words, the DODE problem can be summarized
as the problem of finding the dynamic OD matrices that, when assigned to the
traffic network as vehicles traveling between origins and destinations,
reproduce the traffic measurements captured by sensors in the network. In
addition to this, the dynamic OD matrices should be close to some hypothesized
ones that can be determined, e.g., by the number of people living in the zones
that are considered as origins and destinations.

Solving bi-level optimization problems is quite challenging, as the inner
problem constrains the outer problem, i.e., only an optimal solution to the
inner problem is a feasible solution to the outer problem. To overcome this
challenge, it is crucial to design solution approaches that are efficient and
able to produce good and reliable solutions in a reasonable amount of time. As
such, much of the recent research on the DODE problem tries to include and
leverage additional input data or problem-specific knowledge to improve the
efficiency of the applied algorithms.

The DODE problem is solved frequently for a certain geographical area
of study and traffic network whenever a more up-to-date demand estimate
is needed. If the general structure of the traffic network and the possible OD
pairs within the area of study remain the same over a reasonable amount of time,
the input data in the form of recent traffic observations is the only component
of the DODE problem that changes from one optimization task to another. The
solution approaches that apply sequential optimization algorithms to solve the
DODE problem imply performing many DTAs in each separate optimization task,
essentially starting from scratch every time a new demand estimate is
needed. In a single optimization task, vast computational resources are spent on
performing the DTAs, and each DTA produces valuable information. This
information is usually discarded and not used in subsequent optimization tasks,
primarily due to the inherent sequential nature of some of the applied solution
approaches. In this context, an ML approach becomes appealing as it allows
leveraging data produced by many DTAs.

Therefore, we study an approach to bypass the computationally costly and
time-consuming computations performed in the solution of the inner optimization
problem by learning a model for the input-output relationship established by
solving this inner problem.

To benchmark the ML approaches, we compare them combined with classical
gradient-based approaches. In addition, we provide an analysis of the advantages
and disadvantages of these approaches.

Our work provides the following contributions: We design a possible way to use
ML in the solution of the DODE problem by learning a model for the inner problem
and reducing its computational cost. Further, we perform a fair computational
comparison under a controlled scenario, carefully designed and synthetically
generated, of two state-of-the-art methods and our new ML approach.  

Finally, we provide an analysis of the advantages and disadvantages of the
different approaches and indications about the contexts in which they could be
performing the best.

All methods are here reimplemented and share the same subroutines. We use the
well-established, open-source traffic simulator SUMO \citep{sumo2018} to perform
the DTA in the inner optimization problem. As such, the DTA performed by SUMO
will be regarded as a black-box function that assigns to a given input dynamic
OD matrix a set of path flows in the network according to a certain equilibrium
criterion called a Stochastic User Equilibrium (SUE). This equilibrium criterion
is described in more detail in Section~\ref{sec:problem_statement}. 

The rest of the chapter is outlined as follows. First,
Section~\ref{sec:literature_review} gives a general overview of the existing
literature on the DODE problem. Then, Section~\ref{sec:problem_statement}
introduces notation and the traditional (continuous) bi-level optimization
formulation of the DODE problem that we will use in the subsequent sections. In
Section~\ref{sec:solution_approaches}, we describe the approaches considered in
the study, and in Section~\ref{sec:experimental_setup}, the irregular grid
network and the synthetic problem data that we used for our experiments. Next,
the results are presented in Section~\ref{sec:results}. The best method turns
out to be a gradient-based approach with assignment matrix modeling, while the
ML approach shows inferior performance. Finally, we conclude with a discussion
of the results in Section~\ref{sec:conclusion}.

\section{Literature Review} \label{sec:literature_review}

In the classical static OD estimation problem, static OD trips are considered.
The problem is to estimate the demand in a single time interval (a single
OD matrix is estimated). On the other hand, time-dependent OD trips are
considered in the DODE problem, and the demand can vary in the different
time intervals (i.e., time-dependent OD matrices are estimated).

The static OD matrix estimation problem has a long history, and a large body of
literature on this problem is available. The ideas and insights that have been
gained through studying the static problem are also extended and applied to the
dynamic problem. For an overview of the different classical solution approaches
to the static problem, see \cite{abr1998, mar2009, ber2009}. On the other hand,
the most recent and extensive literature review on the dynamic problem is
provided in \cite{omr2012, dju2014}. In contrast, a literature review on much
earlier work is provided in \cite{bal2006}. For a more general overview of the
DODE problem and the different modeling components and decisions that have to be
considered when solving this problem, we refer the reader to \cite{lin2003,
ber2009}. Finally, an overview of related and relevant problems is provided in
\cite{ant2016}.

Most of the solution approaches to the DODE problem described in the literature
can be divided into two broad classes: DTA-based and non-DTA-based approaches. A
DTA-based approach is usually used in congested networks where the path choices
of the users of the network and the road-use patterns are of interest. DTA-based
approaches provide a sound way to consider behavioral factors that affect travel
decisions but are more computationally demanding than non-DTA-based approaches.
Contrary to this, in a non-DTA-based approach, vehicles are assumed to take the
shortest path, which is usually also the one that minimizes the travel time. In
this case, if the users of the traffic network are assumed to be able to take
the shortest path, then to solve the DODE problem, a non-DTA-based approach is
usually sufficient. Therefore, we focus on the literature concerned with
DTA-based solution approaches with this distinction in mind.

Usually, solution approaches to the outer problem of the bi-level formulation of
the DODE problem consider the problem of minimizing the discrepancy between the
quantities determined by a DTA solution and the corresponding observed
quantities given as input. The problem is treated as a continuous optimization
problem, and general-purpose optimization algorithms from local and global
optimization are used. The main advantage of general-purpose algorithms is that
they do not require exact knowledge of the functional relationship of the
variables to be estimated. The drawback is that these algorithms are usually
inefficient in their original form if the different algorithm components and
parameters are not adapted and tuned to the problem at hand. The general-purpose
algorithms applied to the DODE problem can be categorized as either
\emph{gradient-based} or \emph{derivative-free} optimization algorithms.
Gradient-based algorithms are iterative procedures that adjust estimates based
on the information provided by the gradient and possibly higher-order
derivatives. This is not the case for derivative-free algorithms (e.g.,
simulated annealing, evolutionary algorithms, and sampling-based algorithms,
such as the Nelder-Mead algorithm) that usually rely solely on Objective
Function (OF) evaluations.

We can further identify \emph{assignment matrix-based} or \emph{assignment
matrix-free} approaches within the class of gradient-based algorithms. The
distinction can be made based on whether an Assignment Matrix (AM) is used as a
part of the applied optimization algorithm or not. An AM summarizes the result
of a DTA and consists of elements that describe the utilization of the roads in
the traffic network with respect to the number of users that travel between
origins and destinations within a certain time period. This means that the AM
varies as a function of the travel demands. The solution approaches to the DODE
problem that uses gradients and AMs exploit that it can be possible to
analytically derive the exact gradient from a functional relationship between
the variables through the AM. On the other hand, a gradient-based but AM-free
approach primarily uses finite difference approximations of the gradient, which
are established through OF evaluations.

The most recent studies that adopt a gradient and AM-based approach are given in
\cite{ber2009, tol2013, fred2013, dju2017, sha2017, mas2018}. These studies
primarily focus on improving the efficiency of the applied algorithms by better
modeling the functional relationship between the variables that are to be
estimated. A common theme for the gradient and AM-based approaches is that count
observations are the only type of traffic observations usually given as input to
the optimization problem. The modeling of the relationship between the variables
to be estimated can be done through the AM. However, using this type of approach
limits the type of observations that can be accommodated in the optimization
problem. For example, the functional relationship between demands and
count observations can easily be established if a proportional relationship
between these variables is assumed (see Eq.\ref{eqn:assignmat_relation}).
However, for other types of observations, it can be harder and more cumbersome
to determine appropriate relationships between the variables to be estimated.

Gradient and AM-free approaches have also gained considerable attention in the
literature. Several researchers have especially been interested in studying and
applying different variations of the Stochastic Perturbation Simultaneous
Approximation (SPSA) algorithm by Spall \citep{spall1998}, tailoring it to the
DODE problem. This algorithm only relies on OF evaluations with no need for an
explicit characterization of the relationship between the variables to be
estimated. In each iteration of the SPSA algorithm, OF evaluations are used to
approximate the gradient. This approach makes it possible to include other types
of traffic observations besides count observations, and it has therefore been
the predominant approach in these cases \citep{vaz2009, cip2011, tym2015,
ant2015, lu2015, car2017}. Just like most gradient-based optimization
algorithms, the performance of the SPSA algorithm is sensitive to (i) the tuning
of algorithm parameters, (ii) the possibly very different magnitudes of the
variables to be estimated, and (iii) the OF shape. Different enhancements to the
SPSA algorithm in the mentioned references thus have also been focused on
components that address these points to improve the stability and robustness of
the algorithm when applied to the DODE problem.

A \emph{surrogate model} can be estimated based on empirical data if a model for
the OF is not available through theoretical argumentation. Gradient-based or 
derivative-free methods can then be applied using the  surrogate model.
In \cite{bal2006}, a model-based and derivative-free method is compared against
the SPSA algorithm. Response surface techniques are used to model the OF
by low order polynomials fitted locally to the values in correspondence of
sample points of the search space. The author applies a Stable Noisy
Optimization by Branch and Fit (SNOBFIT) algorithm that uses a derivative-free
Box-Complex algorithm. The Box-Complex algorithm is a direct search and
sampling-based method that extends the well-known Nelder-Mead algorithm. The
Box-Complex and SNOBFIT algorithm performed worse in efficacy and scalability
than the SPSA algorithm.

Generally, the derivative-free solution algorithms require a high number of OF
evaluations to obtain results comparable to the gradient-based solution
algorithms \citep{char2017}.  

A promising line of research is provided in \cite{oso2019}, which describes a
surrogate model-based approach that relies on fitting simplified analytical
models to the underlying simulation and DTA model resulting in a system of
equations that can be solved efficiently by standard system of equations
solvers. The approach thus obtains good computational results under tight
computational budgets. The surrogate model-based approach is compared against
the SPSA and a derivative-free pattern search algorithm. The approach performs
well, while the pattern search algorithm and the SPSA algorithm perform poorly
in comparison.

Our work contributes to this thread of research by studying ML approaches to
determine surrogate models as an alternative to response surface techniques
\cite{bal2006} and the surrogate model-based approach proposed in
\cite{oso2019}.

Fewer studies can be found where model-free (and hence derivative-free)
algorithms have been applied to the DODE problem. These algorithms only rely on
OF evaluations and exhibit the advantage that they can be
distributed and parallelized, as they allow independent OF
evaluations. Evolutionary algorithms have been the most popular among
derivative-free optimization algorithms. Notably, in \cite{kat2006}, an
evolutionary algorithm was applied in a setting with distributed and parallel
computing. Other studies that apply evolutionary algorithms can be found in
\cite{vaz2009,omr2014}, while \cite{tse2007} also use an evolutionary algorithm,
but in a non-DTA based approach. However, in \cite{vaz2009}, the evolutionary
algorithm was shown to perform worse than the SPSA algorithm.

Two variants of DODE problems are usually addressed in the literature. Their
offline or online setting characterizes them. In an offline setting, historical
data is used to estimate the network-wide demand, such that long-term
predictions of future traffic conditions can be made. This is especially
beneficial to decision-makers. It enables them to evaluate different planning
scenarios and make informed decisions when proposed changes to infrastructure
and facilities are made. Offline demand estimation procedures are usually
applied in network-wide studies, so the computational cost can be considerably
high in this setting. In an online setting, real-time traffic data is used and
historical data to estimate the current level of demand such that short-term
predictions of future traffic conditions can be made. This primarily benefits
real-time traffic control and path-guidance systems. Online demand estimation
procedures are usually applied on smaller traffic networks or locally, e.g., at
intersections, where traffic patterns are identified, such that adaptive traffic
control strategies can be used. For references on work that studies the online
DODE problem, we refer the reader to \cite{ash2000, ant2009, omr2012}. In this
work, we focus on the offline, DTA-based, bi-level formulation of the DODE
problem.

\section{Notation and Problem Statement} \label{sec:problem_statement}

A traffic network can be defined as a directed graph $G = \left(N, A\right)$
consisting of a set of nodes $N$ and a set of arcs $A$. Each node in the network
represents a junction or an origin/destination point, while arcs represent
roads. The arcs are directed, which means that the traffic between two nodes can
be uni-directional and one-way roads are possible. The arcs in the network have
additional attributes that generally indicate how much and how well the roads in
the network can handle traffic. Some examples of essential arc attributes
are the number of lanes, the length $\left(\textnormal{kilometer}\right)$ and a
free-flow speed
$\left(\nicefrac{\textnormal{kilometer}}{\textnormal{hour}}\right)$ that is the
maximum allowed speed a vehicle can travel with on an uncongested road.
Additionally, geographical information can be added to the arcs and nodes to
indicate their physical location and shape.

A subset of the nodes $N$ in a traffic network can be identified as possible
origins and/or destinations. More precisely, we let $O \subseteq N$ be a set of
origins and $D \subseteq N$ a set of destinations, where it is usually the case
that $D \cap O \ne \emptyset$, i.e., it is possible for a node to be an origin
and a destination simultaneously. A \emph{trip} is the movement of a vehicle
from one location to another. More precisely, a trip departs from an origin $i
\in O$ and terminates at a destination $j \in D$. In this case, we associate the
trip with an OD $w = (i, j) \in W  = O \times D$. Here $W$ is the set of
all OD pairs with size $\left|W\right| = m$, and it is assumed
that no trip departs and arrives in the same location, i.e., the OD pairs are
defined such that $i \neq j$ for all $w = (i, j) \in W$. Finally, a subset of
the arcs in the network $G$ is assumed to be equipped with sensors defined by
the set $Q = \left\{1, \hdots, n_Q\right\} \subseteq A$.

To define the quantities of interest in the DODE problem, we introduce an
analysis period $T = \left[0, t_{\textnormal{end}}\right)$ discretized into a
set of $n_S$ subintervals $\left\{[0, t_1), [t_1, t_2), \dots, [t_{n_S - 1},
t_{n_S} = t_{\textnormal{end}})\right\}$ of equal duration. For the sake of
convenience, we identify these intervals by the indices in the set $S =
\left\{1,\hdots, n_S\right\}$. These indices define the time intervals in which
we want to estimate the demands for each of the $m$ OD pairs. The estimation of
the demands is then based on the arc count observations made within these time
intervals.

The quantities of interest are usually referred to as being arranged in
matrices. Here, however, we will arrange these quantities in vectors. In other
words, we flatten or vectorize the matrices, i.e., given a matrix $\mathbf{A}
\in \mathbb{R}_{+}^{p_1 \times p_2}$, we concatenate consecutive rows of the
matrix in a column vector ($\mathbb{R}^{p_1 \times p_2} \rightarrow
\mathbb{R}^{p_1 \cdot p_2}$):
\begin{align}
    \textnormal{vec}(\mathbf{A}) = \mathbf{a} = \left[a_{11},\dots, a_{p_11},\dots, a_{12}, \dots, a_{p_1 2}, \dots, a_{p_1 p_2}\right]^\top \in \mathbb{R}^{p_1 \cdot p_2 }.
\end{align} 

The quantities of interest can then be described in terms of the following
vectors:

\begin{itemize}
    \item $\mathbf{x}, \tilde{\mathbf{x}} \in \mathbb{R}_{+}^{m \cdot n_S}$ are
    the vectors of \emph{estimated demands} and \emph{seed demands},
    respectively. Each element of these vectors ${x}_{ws}$ or $\tilde{x}_{ws}$
    is associated with an OD pair $w \in W$ and a time interval $s \in S$. The
    vector $\mathbf{x}$ of estimated demands contains the estimate of the number
    of trips made for each OD pair in each of the time intervals. The vector
    $\tilde{\mathbf{x}}$ of seed demands defines the prior knowledge of the
    number of trips made for each OD pair in each time interval. This prior
    knowledge is usually assumed to have been obtained from a previous demand
    study, e.g., from a population survey. Note that we model discrete values as
    real numbers. We thus define the continuous relaxation of a discrete
    optimization problem. 
    
    \item $\mathbf{x}^{\text{Lower}}, \mathbf{x}^{\text{Upper}} \in
    \mathbb{R}_{+}^{m \cdot n_S}$ are the vectors of lower and upper bounds on
    the estimated demands, respectively. These upper and lower bounds define the
    \emph{search space}.
    \item $\mathbf{c}, \hat{\mathbf{c}} \in \mathbb{Z}_{+}^{n_Q \cdot n_S}$ are
    the vectors of estimated and observed arc counts, respectively. Each element
    of these vectors $c_{qs}$ or $\hat{c}_{qs}$ is associated with a sensor $q
    \in Q$ and a time interval $s \in S$.
\end{itemize}

The outer optimization problem of the bi-level formulation of the DODE problem
can now be defined as:
\begin{align}
    \min\;\;& F(\mathbf{x}, \tilde{\mathbf{x}}, \mathbf{c}(\mathbf{x}), \hat{\mathbf{c}})
    \label{eqn:1_estimation_problem} \\
    \textrm{subject to } \;\; & \mathbf{x}^{\text{Lower}} \leq \mathbf{x} \leq \mathbf{x}^{\text{Upper}}
    \label{eqn:2_estimation_problem}
\end{align}

The solution of this problem, $\mathbf{x}^*$, is the demand estimate that
results from minimizing the OF $F$, which we define as the
weighted sum of the measures of discrepancy between the estimated quantities and
their corresponding observed or a priori values:
\begin{align}
    F(\mathbf{x}, \tilde{\mathbf{x}}, \mathbf{c}(\mathbf{x}), \hat{\mathbf{c}}) = \omega_1 \cdot f^{(1)}(\mathbf{x}, \tilde{\mathbf{x}}) + \omega_2 \cdot f^{(2)}(\mathbf{c}(\mathbf{x}), \hat{\mathbf{c}}), \quad \omega_1,\omega_2\in \mathbb{R}_+.
    \label{eqn:objective_function}
\end{align}

In our specific case, the discrepancy $f^{(1)}$ is measured between the
estimated demands $\mathbf{x}$ and a priori known demands
$\tilde{\mathbf{x}}$ and the discrepancy $f^{(2)}$ between the estimated arc
counts $\mathbf{c}(\mathbf{x})$ and the observed arc counts $\hat{\mathbf{c}}$.
Note that in a more general setting, the OF may consist of
several additional terms $f^{(3)}, f^{(4)}, \dots$ with respective weights
$\omega_3, \omega_4, \dots$ that consider additional available information, such
as arc speed, density, travel times, turning ratios, etc., that can be included
to improve the guidance of the search process.

An inner optimization problem must be solved to evaluate the OF in
Eq.~\eqref{eqn:1_estimation_problem} for a proposed vector $\mathbf{x}$ of
demands to find a corresponding vector $\mathbf{c}(\mathbf{x})$ of arc counts.
The inner optimization problem takes the form of a DTA problem. The task is to
determine the routing of users between origins and destinations according to a
particular assignment principle. The users of the traffic network are assumed to
make decisions per the criterion specified by the assignment principle. Two
examples of assignment principles that are widely used and studied within the
area of transportation research are (i) the User Equilibrium (UE) principle
(Wardrop's first principle \citep{fri2016}), which states that each user makes
decisions that minimize their own individual travel time, and (ii) the System
Optimum (SO) principle (Wardrop's second principle \citep{fri2016}), which
states that users make joint decisions to minimize the total system travel time.

DTA models adopt an assignment principle and incorporate several different
travel choice components to reflect the travel choices a user of an actual
traffic network might face when wanting to travel between an origin and a
destination. Travel choices that are possible to model and include in a DTA
model are, among others: path choice, departure time choice, mode choice, and
destination choice, to name a few. However, in the context of the DODE problem,
a simple DTA model that only incorporates path and departure time choice
components is usually used, meaning path and departure time choices are
endogenous to the DTA model, while other travel choice components are exogenous
or fixed.

To obtain accurate travel times between origins and destinations, a DTA model
relies on an underlying traffic flow model, i.e., detailed traffic flow models
are used by DTA models to explicitly propagate vehicles from one arc to another
while considering space constraints. In other words, queuing and congestion are
modeled, and the effects of these phenomena are reflected in the travel time.
For a good introduction and a general overview of the DTA problem, we refer the
reader to \cite{peeta2001}, as we only give a brief description here.

Two traffic flow models are implemented in SUMO and can be used in conjunction
with a DTA: a \emph{mesoscopic} and a \emph{microscopic}. These two types of
traffic flow models differ in the level of detail in which they model traffic
flow dynamics. In the mesoscopic model, the arcs in a traffic network are
modeled as discrete queues through which vehicles are propagated. Further, a
coarser model for intersections and lane-changing is used \citep{sumo_meso}. On
the other hand, in the microscopic model, detailed quantities such as position,
speed, acceleration, and braking distance are considered in the propagation of
vehicles on an arc.

The mesoscopic model is computationally cheaper to solve and can provide the
necessary data to a sufficient degree of detail needed in the estimation
problem. Therefore, we choose that model. 

A DTA can be obtained heuristically through an iterative simulation process
using the mesoscopic traffic flow model implemented in SUMO. At the end of the
iterative process, a stationary distribution of the path choice decisions of the
users of the traffic network respects an assignment principle. This assignment
principle is the Stochastic User Equilibrium (SUE) in SUMO. The SUE assignment
principle extends the UE assignment principle, where stochastic elements have
been incorporated into the DTA model. This form of UE assignment is regarded as
the more realistic, as uncertainty is incorporated in the assignment model to
take into account the uncertainty of the users' knowledge about network
conditions. Contrary to the UE assignment principle, which is defined in terms
of \emph{actual} travel times, the SUE assignment principle is instead defined
in terms of \emph{perceived} travel times.

The iterative method used by SUMO to determine a DTA that respects the SUE
assignment principle is described in \cite{gaw1998}. This method determines the
probabilities of choosing between path alternatives for each user who wants to
travel between an origin and destination at a certain time. In brief, the path
choice probabilities are determined based on (i) the arc travel times
experienced in the previous iteration, (ii) the sum of arc travel times along
different least-cost paths (these paths constitute a set of alternatives to a
user), and (iii) the previous probabilities of choosing the paths. These
quantities are used to obtain new estimates of the path choice probabilities at
each iteration. Finally, we note that the least-cost paths are computed
in parallel in SUMO by a time-dependent version of Dijkstra's algorithm, 
while a single traffic simulation is single-threaded and sequential.

Note that we modeled demands $\mathbf{x}$ as continuous variables while
the mesoscopic simulator implemented by SUMO solves a discrete optimization
problem. When the demands are given to SUMO as input, they are rounded to
the nearest integer values. In all other cases, the demands are handled
as being continuous, which allows us to apply continuous optimization algorithms
in the outer-level problem.

For our purposes here, the DTA performed by SUMO can conveniently be defined by
the function:
\begin{align} \label{eqn:sumo_dta_map}
    \Gamma: \mathbb{R}_{+}^{m \cdot n_S} \rightarrow \mathbb{Z}_{+}^{n_Q \cdot n_S}.
\end{align}
This function maps demands $\mathbf{x}$ to arc counts through a DTA, i.e.,
$\mathbf{x} \mapsto \Gamma(\mathbf{\mathbf{x}}) = \mathbf{c}(\mathbf{x})$. The
arc counts $\mathbf{c}(\mathbf{x})$, together with the corresponding input
demands $\mathbf{x}$, can then be used in the OF in
Eq.~\eqref{eqn:objective_function} to determine the quality of the current 
demand estimate.

Eq.~\eqref{eqn:2_estimation_problem} defines upper and lower bounds on the
demands. These are usually supplied to avoid a situation where an unrealistic
high amount of demand is loaded on the traffic network, exceeding the network's
capacity. In a realistic setting, the bounds on the demands are usually
determined based on population survey data and experimental results from
previous studies.

In the remaining part of this work, we use the shorthand notation
$F(\mathbf{x})$ for the OF defined in Eq.~\eqref{eqn:objective_function}, as  
$\tilde{\mathbf{x}}$ and $\hat{\mathbf{c}}$ are given and
$\mathbf{c}(\mathbf{x})$ is derived from $\mathbf{x}$.

\section{Solution Approaches} \label{sec:solution_approaches}

Different definitions of the discrepancy between observed and estimated
quantities in the OF $F$ are possible. Drawing on the knowledge obtained in
\cite{cip2011} and \cite{ant2016}, we use an OF that penalizes large errors
between observed and estimated quantities and that weights the different terms
evenly, i.e., we take into account that the different quantities $f^{(1)}$ and
$f^{(2)}$ that enter into the OF can have different magnitudes and hence
normalize their values. Thus, we define the terms $f^{(1)}$ and $f^{(2)}$ used
in the OF $F$ in the following functional form:

\begin{align}
    f^{(1)}(\mathbf{x}, \tilde{\mathbf{x}}) =
        \frac{
            \sqrt{
            \sum\limits_{ \substack{w \in W \\ s \in S}  } \left(x_{ws} - \tilde{x}_{ws}\right)^2
            }
        }
        {
            \sqrt{\sum\limits_{  \substack{w \in W \\ s \in S}  } \tilde{x}_{ws}^2}
        }
\;\; \textnormal{and}\;\;
f^{(2)}(\mathbf{c}(\mathbf{x}), \hat{\mathbf{c}}) =
    \frac{
        \sqrt{
            \sum\limits_{ \substack{q \in Q \\ s \in S}  } \left(c_{qs} - \hat{c}_{qs}\right)^2
        }
    }
    {
        \sqrt{\sum\limits_{  \substack{q \in Q \\ s \in S}  } \hat{c}_{qs}^2}
    }.
\label{eqn:objective_function_terms}
\end{align}

To further limit the search space, we introduce generation constraints
\citep{ant2016} that set an upper bound on the number of outbound trips from an
origin.  we let $o_i$ denote the maximum number of vehicles that can leave
origin $i \in O$ and enforce that:
\begin{align}
    \sum_{\substack{j \in D, \; w = (i, j) \in W \\ s \in S}} x_{ws} \leq o_i \quad \forall i \in O.
    \label{eqn:generation_constraints}
\end{align}

This type of constraint can easily be incorporated into the problem formulation
as it is typically the case that data for the values of $o_i$, $i\in O$ are
readily available and highly reliable, e.g., from survey data or population
density statistics. Other possible constraints are also mentioned in
\cite{ant2016}, such as trip distribution constraints or attraction constraints
that, contrary to the generation constraints, limit the number of trips that can
arrive at a certain destination. However, data pertaining to these constraints
might be harder to provide and will thus not be considered in this study.

The algorithms that follow are iterative algorithms. We define an iteration
index $\tau \in \left\{0, \hdots, \tau_{\textnormal{max}}\right\}$ where
$\tau_{\textnormal{max}}$ is the maximum number of iterations of an algorithm.
We thus let $\mathbf{x}_{\tau}$ define the estimate at iteration $\tau$.

\subsection{Gradient-based Approaches} \label{subsec:gradient_based_approach}

Gradient-based algorithms for continuous optimization problems
take the general form of an iterative procedure where an incumbent estimate at
iteration $\tau$ is adjusted in such a way that it is moved towards a minimum to
yield the new estimate at iteration $\tau + 1$:
\begin{align}
    \mathbf{x}_{\tau+1} = \mathbf{x}_{\tau} + \eta_{\tau} \cdot \mathbf{d}_{\tau}.
    \label{eq:adjustment_step}
\end{align}
Here $\eta_{\tau}$ is a scalar value that we will refer to as the
\emph{step rate}, and $\mathbf{d}_{\tau}$ is the \emph{descent
  direction}. The descent direction is determined by the gradient
through $\mathbf{d}_{\tau}=-\mathbf{g}_{\tau}$. Different gradient-based
optimization algorithms differ in how the descent direction and the step
rate are chosen. For the DODE problem, it is possible to distinguish
two ways to derive the gradient: AM-based and AM-free.

\subsubsection{An Assignment Matrix-based Approach} \label{sec:assignmat_based_apporach}

Through the use of an AM the gradient of the OF
can be calculated analytically. The (vectorized) AM
$\mathbf{p}_{\tau} = \textnormal{vec}(\mathbf{P}_\tau(\mathbf{x}_{\tau}))$ at an
iteration $\tau$ consists of entries $p_{qrws} \in \left[0, 1\right]$ that
describe the proportion of trips $x_{ws}$ for OD pair $w \in W$ in time interval
$s \in S$ that go through arc $q \in Q$ in time interval $r \in S, \; r \geq s$.
In other words, depending on a number of factors such as the travel time between
origins and destinations and the length of the time intervals $s \in S$, the
traffic present on an arc $a \in A$ may originate from OD trips from several
different OD pairs from the current time interval or several previous time
intervals. Using this knowledge the OD trips can be related to the arc counts
through the linear relation:
\begin{align}
      c_{qr} = \sum_{\substack{w \in W \\ s \in S \\ s \leq r}} p_{qrws} \cdot x_{ws}, \qquad \forall q \in Q, r \in S.
      \label{eqn:assignmat_relation}
\end{align}
Through the use of the AM it is possible to derive an analytical
expression for the gradient of the OF $F$ at every iteration:
\begin{align}
    \frac{\partial F}{\partial \mathbf{x}} &= \omega_1 \cdot \mathbf{g}^{(1)} + \omega_2 \cdot \mathbf{g}^{(2)},
    \label{eqn:original_descent_direction}
\end{align}
where the vectors $\mathbf{g}^{(1)}$ and $\mathbf{g}^{(2)}$ are the gradients of
the two terms  $f^{(1)}$ and  $f^{(2)}$, respectively, whose elements can be
derived to be (see \ref{sec:grad_derivation}):

\begin{align}
    g^{(1)}_{ws} = \frac{ \partial f^{(1)}  }{  \partial x_{ws}  }  
    =
        \frac{
            x_{ws} - \tilde{x}_{ws}
        }
        {
            \sqrt{
                \sum\limits_{  \substack{w' \in W \\ s' \in S  }  } \left(x_{w's'} - \tilde{x}_{w's'}\right)^2
            } \cdot \sqrt{\sum\limits_{  \substack{w' \in W \\ s' \in S}  } \tilde{x}_{w's'}^2}
        }, \quad \forall w \in W, s \in S,
    \label{eqn:first_grad}\\[2em]
    g^{(2)}_{ws} = \frac{ \partial f^{(2)}  }{  \partial x_{ws}  } 
    = \frac{
        \sum\limits_{  \substack{q' \in Q \\ r' \in S\\ r'\geq s}  }
            \left(
                c_{q'r'}
                - \hat{c}_{q'r'}
            \right) \cdot p_{q'r'ws}
    }
    {
        \sqrt{
             \sum\limits_{  \substack{q' \in Q \\ r' \in S  }  } \left(c_{q'r'} - \hat{c}_{q'r'}\right)^2
        } \cdot \sqrt{\sum\limits_{  \substack{q' \in Q \\ r' \in S}  } \hat{c}_{q'r'}}
    }, \quad \forall w \in W, s \in S.
    \label{eqn:second_grad}
\end{align}

The descent direction in Eq.~\eqref{eq:adjustment_step} is then simply the
negative of the expression in Eq.~\eqref{eqn:original_descent_direction}: 

\begin{align}
    \mathbf{d}_{\tau} =
        - \left(
            \omega_1 \cdot \mathbf{g}^{(1)}_\tau + \omega_2 \cdot \mathbf{g}^{(2)}_\tau
        \right).
        \label{eqn:descent_direction_assignmat_based}
\end{align}

\subsubsection{An Assignment Matrix-free Approach}
\label{sec:assignment_matrix_free}

As mentioned in Section~\ref{sec:literature_review}, most of the studied
gradient-based and AM-free approaches for the DODE problem use the SPSA
algorithm. This algorithm uses the negative of the gradient as a way to obtain a
descent direction in a similar way as the previous approach. As no AM is needed
in the derivation of the gradient of the first OF term $f^{(1)}$, we can simply
use the exact gradient defined in Eq.~\eqref{eqn:first_grad}. However, the
gradient of the second OF term $f^{(2)}$ will have to be approximated because
now no AM is given and no functional relationship between the variables to be
estimated is assumed. 

A one-sided finite difference approximation of the gradient is usually used:
\begin{align}
    \tilde{\mathbf{g}}^{(2)}_{\tau} = \frac{f^{(2)}(\mathbf{c}(\mathbf{x}_\tau), \hat{\mathbf{c}}) - f^{(2)}(\mathbf{c}^{-}(\mathbf{x}_\tau), \hat{\mathbf{c}})}{b_\tau}
    \cdot
    \overline{\mathbf{\Delta}}_{\tau}
    , \quad \tau \in \left\{0, \dots, \tau_{\textnormal{max}} \right\}.
    \label{eqn:grad_approx_1}
\end{align}
Here, $\mathbf{c}(\mathbf{x}_\tau) = \Gamma(\mathbf{x}_{\tau})$, $
\mathbf{c}^{-}(\mathbf{x}_\tau) = \Gamma(\mathbf{x}_{\tau} - b_{\tau} \cdot
\mathbf{\Delta}_\tau)$ and $\mathbf{\Delta}_\tau \; \in \left\{-1, 1\right\}^{m
\cdot n_S}$ is a random perturbation vector, where the entries $\Delta_{ws}, w
\in W, s \in S$ are independent Bernoulli random variables taking the value $\pm
1$ with probability $\nicefrac{1}{2}$. Furthermore, the vector
$\overline{\mathbf{\Delta}}_\tau$ consists of entries
$\nicefrac{1}{\Delta_{ws}}$ and $b_\tau = \nicefrac{b}{\tau^\gamma}$ is defined
with respect to the parameters $b > 0$ and $\gamma > 0$. These parameters and
the Bernoulli distribution are usually chosen because they satisfy certain
conditions \citep{spall1992} that ensure the algorithm to be able to converge to
a minimum, asymptotically. To ensure compatibility with SUMO that only accepts
integer inputs, we keep the perturbation parameter fixed $b_\tau = 1$ throughout
all iterations, which is equivalent to perturbing each element in the demand
vector by $\pm 1$ trip.

\subsubsection{Determining Step Lengths} \label{sec:determine_step_length}

The step length in the descent direction can be chosen appropriately  
by solving  for $\eta_\tau$  a line search problem in each iteration $\tau$:
\begin{align}
    \eta_\tau = \arg \min_{\alpha \geq 0} F\left(\mathbf{x}_{\tau} + \alpha \cdot \mathbf{d}_{\tau}\right).
    \label{eqn:line_search}
\end{align}

The goal of the line search is to find the step rate $\alpha$ that minimizes the
OF in the descent direction given by $\mathbf{d}_\tau$. To solve this
sub-problem we use a multiple linear regression approach where a 2nd-degree
polynomial is fitted to sample points evaluated in the descent direction. More
specifically, we first define a single-variable function of the step rate
$\alpha$, as follows:
\begin{align}
    y(\alpha) = F\left(\mathbf{x}_{\tau} + \alpha \cdot \mathbf{d}_\tau\right).
    \label{eqn:uni_var}
\end{align}
Then, we determine:
\begin{enumerate}
    \label{list:ub}
    \item An upper bound $\ell$ on $\alpha$ imposed by the generation constraints in
      Eq.~\eqref{eqn:generation_constraints}. An algorithmic sketch to
      carry out this task is given in
      Algorithm~\ref{alg:determine_max_step_length}.
    \item A number of equally spaced points $\alpha_{1}, \hdots, \alpha_n, \; n
    \geq 2$ to evaluate in the interval $\left[\nicefrac{\ell}{n}, \ell\right]$.
\end{enumerate}
Given these parameters we can obtain sample points with corresponding responses:
\[(\alpha_0, y(\alpha_0)), \hdots, (\alpha_n, y(\alpha_n)).\] The point
$\left(\alpha_0 = 0,\; y(\alpha_{0}) = F\left(\mathbf{x}_{\tau}\right)\right)$
is already known. In this case, by specifying two additional sample points we
are able to uniquely fit a 2nd-degree polynomial. More precisely, an estimator
for the function $y(\alpha)$ in Eq.~\eqref{eqn:uni_var} is:
\begin{align}
    \hat{y}(\alpha) = \beta_0 + \beta_1 \cdot \alpha + \beta_2\cdot \alpha^2
\end{align}
where the coefficients $\beta_0, \hdots, \beta_2$ are determined by \emph{least
squares}, given the sample points and corresponding responses. The value
$\alpha_{\textnormal{min}}$ that minimizes this polynomial is given by
$\alpha_{\textnormal{min}} = \frac{-\beta_{1}}{2 \cdot \beta_{2}}$, if $\beta_2
> 0$. Otherwise, $\alpha_{\textnormal{min}}$ is taken to be the step rate that
produced the sample point with the smallest response.

\IncMargin{1em}
\begin{algorithm}[tb]
    \small
    \newcommand\mycommfont[1]{\ttfamily{#1}}
    \SetCommentSty{mycommfont}
    
    \SetAlgoLined\DontPrintSemicolon
    \SetKwFunction{proc}{DetermineMaxStepRate}
    \SetKwProg{myprocedure}{Procedure}{}{}
    \myprocedure{\proc{$\mathbf{x}_{\tau}, \,\mathbf{x}_{l}, \,\mathbf{x}_{u},
    \,\mathbf{d}_{\tau}, \,\ell_{\Delta} = 10^{8}, \,\epsilon_1 = 10^{-8}$}}{
    $\ell \leftarrow 0$\; \While{\texttt{ true }}{ $\ell \leftarrow \ell +
    \ell_{\Delta}$ \;            
            $\breve{\mathbf{x}} \leftarrow \mathbf{x}_{\tau} + \ell \cdot
            \mathbf{d}_{\tau}$\;            
            % Apply upper and lower bounds
            \tcp{Apply upper and lower bounds }
            $\breve{\mathbf{x}} \leftarrow \texttt{maximum(}\mathbf{x}_{l},
            \,\,\breve{\mathbf{x}}\texttt{)}$\; $\breve{\mathbf{x}} \leftarrow
            \texttt{minimum(}\mathbf{x}_{u}, \,\,\breve{\mathbf{x}}\texttt{)}$\;

            % If generation constraints are violated or upper and lower bound
            % constraints are binding
            \tcp{Check if generation constraints are violated or whether }
            \tcp{all upper or lower bound constraints are binding}
            \If{                
                $
                \sum\limits_{\substack{j \in D,\,\, w = (i, j) \in W \\ s \in S}} \breve{x}_{ws} \texttt{ > } o_i \,\, \forall i \in O
                \,\, 
                \vee
                $
                \mbox{}\phantom{\textbf{if}}
                $ 
                % \,\,
                % \left(\breve{\mathbf{x}} == \mathbf{x}_{l} \vee \breve{\mathbf{x}} == \mathbf{x}_{u}\right)
                (\breve{x}_{ws} \texttt{ == } x^{\textnormal{Lower}}_{ws} \vee \breve{x}_{ws} \texttt{ == } x^{\textnormal{Upper}}_{ws}) \,\, \forall w \in W, s\in S
                $                
            }{ \tcp{Backtrack and try again with a smaller step rate }            
                $\ell \leftarrow \ell - \ell_{\Delta}$\; $\ell_{\Delta}
                \leftarrow \nicefrac{\ell_{\Delta}}{2}$\; \If{$\ell_{\Delta}
                \texttt{ < } \epsilon_1$}{ \texttt{break}\; }
            % else
            } } }
    \tcp{Return the largest positive step rate that does not violate any constraints}
    \KwRet $\ell$ \;
    \caption{Procedure for determining the upper bound on the step rate.}
    \label{alg:determine_max_step_length}
\end{algorithm} 
\DecMargin{1em}

\subsubsection{Computational Complexity of the Gradient-based Approaches}

If we look at the computational cost associated with the AM-based approach, we
notice that two OF evaluations and thus DTA problems are solved in each
iteration. First, a single DTA is needed to obtain a (vectorized) AM
$\mathbf{p}_{\tau}$ such that the descent direction and two additional points
(using the relation in Eq.~\eqref{eqn:assignmat_relation}) in the descent
direction can be evaluated, and a 2nd-degree polynomial can be fitted. Using the
2nd-degree polynomial, the best step length in the descent direction can be
determined. Then, one additional DTA is solved to evaluate the new point that
has been found by using the best step length that minimizes the OF in the
descent direction.

On the other hand, the AM-free approach needs to solve up to four DTA problems
in each iteration, i.e., one DTA is required to estimate the gradient, and two
additional DTAs are needed for determining the best step length in the descent
direction. As with the AM-based approach, an additional DTA is then necessary to
evaluate a new estimate found.

\subsection{Machine Learning Approaches} \label{subsec:ml_approach}

The DTA performed by SUMO can be regarded as a black-box function that maps
input demands to network-specific quantities, such as arc counts, speeds,
densities, etc. As mentioned in Section~\ref{sec:problem_statement}, we will
focus on the arc counts. In this case, the black-box function representing the
DTA performed by SUMO can be described by Eq.~\eqref{eqn:sumo_dta_map}, a
vector-valued multivariate function that describes the input-output relationship
of the DTA performed by SUMO with respect to the output quantity we are
interested in modeling. The ML approach's main idea is to learn a model of this
input-output relationship.

A dataset is needed to learn a model of the input-output relationship. This
dataset can be obtained through a \emph{sampling strategy}, which is a technique
for deciding how many and which points should be included in a dataset provided
to a \emph{learning algorithm}. Using the dataset, the learning algorithm is
trained in a supervised fashion, which results in a model $\hat{\Gamma}$ for the
relationship in Eq.~\eqref{eqn:sumo_dta_map}. Subsequently, the DODE problem can
be solved using the obtained model as a surrogate for DTA in the OF.  

\begin{figure*}%[tb!]
    \includegraphics[width=\textwidth,keepaspectratio]{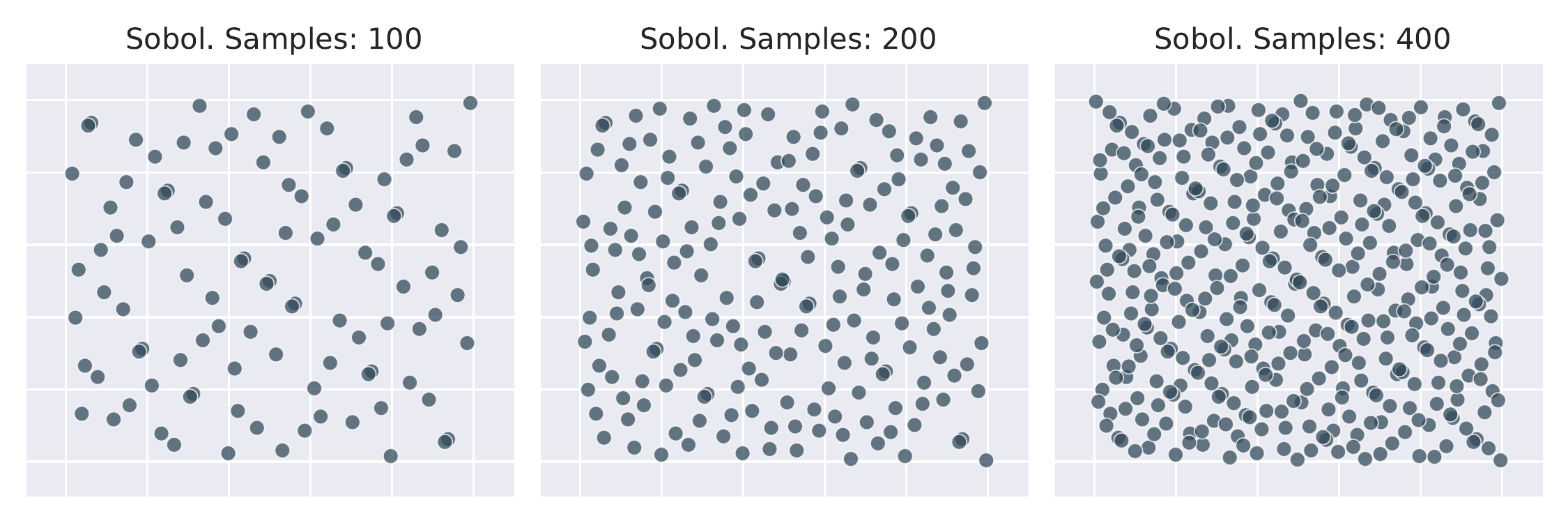}
    \caption{Example of a sampling strategy in the plane that uses a Sobol low-discrepancy sequence.}
    \label{fig:sampling_methods}
\end{figure*}

Overall, the ML approach put forth here consists of three primary
components: a sampling strategy, the application of a ML algorithm
and a final application of an optimization algorithm to obtain an approximate
minimum of the DODE problem that is then to be evaluated and returned as a final
estimate. These steps are described in more detail in the following.

\subsubsection{Sampling Strategy}

The sampling strategy should preferably be chosen such that the amount of
information gained from a limited number of sample points is maximal. Sampling
more points than necessary implies additional computational resources and time.
In other words, for the dataset we want to obtain a set of points that are
sampled evenly in the search space. A way to accomplish this is by using
\emph{optimal experimental designs}, which among others include: Box-Behnken,
central composite, factorial designs, and latin hypercube sampling, to name a
few. These sampling strategies are specialized because they are tailored for
fitting specific parametric models (e.g., linear or multiple linear models).
Moreover, they dictate the exact number of samples needed to end up with a set
of sample points that covers the search space well. Consequently, these
strategies are well-suited for fitting models when prior knowledge of the
model's structure is available, and the sampling budget allows for the number of
samples required by the sampling strategy. If this is not the case, an
alternative is to use a strategy that samples points according to a
\emph{low-discrepancy sequence} (also referred to as a quasi-random sequence).
Examples of low-discrepancy sequences include Halton, Hammersley, Niederreiter,
and Sobol (see \cite{par2016} for more theoretical details on low-discrepancy
sequences).

An advantage of the strategies that generate sample points according to
low-discrepancy sequences is that they can evenly cover the search space,
regardless of the number of samples. It can also be possible to add additional
samples later to fill the search space progressively.

Ultimately, the chosen sampling strategy depends on (i) the dimensionality of
the problem and (ii) the simulation budget, which imposes an upper bound on the
number of points to evaluate, i.e., the number of DTAs to perform. 

In our experiments (see Section~\ref{sec:results}), we will use a dataset with
points sampled according to the numbers in a Sobol sequence (an example is shown
in Figure~\ref{fig:sampling_methods})).

Formally, let $\mathcal{D}_\Gamma = \left\{(\mathbf{x}_{1}, \mathbf{y}_{1}),
\hdots, (\mathbf{x}_{n}, \mathbf{y}_{n}) \right\}$ denote the dataset
where $\mathbf{x}_{i} \in \mathbb{R}_+^{m \cdot n_S}$ is an element of a Sobol
sequence and $\mathbf{y}_{i}= \Gamma(\mathbf{x}_{i}) \in \mathbb{R}_+^{n_Q \cdot
n_S}$ for all $i \in \left\{1, \hdots, n\right\}$. Recall that $\Gamma$
represents the function applied by the SUMO simulator. The learning task is
to determine $\hat\Gamma$ such that it approximates well $\Gamma$ on new data
points.

\subsubsection{Machine Learning Algorithms}

A wide selection of off-the-shelf implementations of general-purpose ML
algorithms can be used to find $\hat\Gamma$. 
% We refer the interested readers to
% the book \emph{The Elements of Statistical Learning} \citep{has2009}. 
The ML algorithms we are interested in studying can be categorized as
nonparametric models. In contrast to parametric models, nonparametric models are
more flexible as no prior assumptions are made on the structure of the model
that is to be trained.

It could be possible to learn the input-output relationship between the demand
and the scalar value of the second OF term in
Eq.~\eqref{eqn:objective_function}, that is:
\begin{align}
    f^{(2)} : \mathbb{R}_+^{n_Q \cdot n_S} \rightarrow \mathbb{R}_+,
    \label{eqn:objective_function_map}
\end{align}
This term of the OF measures the discrepancy between the observed arc counts and
the current estimate of arc counts, which depends indeed on the DTA. On the
contrary, the other term in Eq.~\eqref{eqn:objective_function}, $f^{(1)}$,
measures the discrepancy between the estimated demands and the seed demands and
does not require performing a DTA.

However, preliminary experiments showed that a better approach is trying to
learn the input-output relationship of a vector-valued multivariate function as
in Eq.~\eqref{eqn:sumo_dta_map} and predict the vector of arc counts. Note that
this approach has the advantage of being more robust to a change in the number
of sensors in the network (e.g., due to failure), since in that case the
corresponding predicted arc counts can be ignored in the calculation of
$f^{(2)}$. Learning a single model to predict $p$ target variables can be done
by inherently multi-target ML algorithms, which are able to exploit inter-target
correlations. The inherently multi-target ML algorithms that we will make use of
are \emph{Feed-Forward Networks (FNNs)} and $k$-Nearest Neighbors ($k$-NN). The
algorithms are briefly described in the following and in more detail in
Appendix~\ref{sec:ml_models}.

A \emph{Feed-Forward Network (FNN)} is a type of neural network consisting of
\emph{units} arranged in layers: An input layer, a number of hidden layers and
an output layer. A FNN architecture is thus mainly defined by the number of
layers (the depth of the network $m_d \in \mathbb{N}$ not counting the input
layer) and the number of units in each layer (the width of a layer $l_i \in
\mathbb{N}$, with $i \in \left\{1, \dots, m_d \right\}$), along with the
\emph{activation functions} used in the different layers. In the network
architectures that we consider a \emph{ReLU} ($\phi(\cdot) = \max(0, \cdot)$)
\emph{activation function} is used in the units situated in the hidden layers
and a linear activation function is used in the output layer, since the task of
learning $\Gamma$ is a regression task. Finally, to reduce overfitting and
improve the out-of-sample predictive performance a regularization parameter
$\lambda \in \mathbb{R}_+$ is also set. An appropriate value for this parameter
is usually chosen using out-of-sample \emph{Cross-Validation} (CV) techniques.   

\emph{$k$-NN} is one of the simplest nonparametric models used for regression.
To predict the value $\hat{\mathbf{y}}$ of a new unseen point $\mathbf{x}$, $k$
needs to be specified. Its value specifies the number of points from a training
set $\mathcal{D}$ closest to $\mathbf{x}$ that should be used to calculate the
predicted value $\hat{\mathbf{y}}$. Unlike the FNN and most other ML algorithms,
no model must be determined, rather the training data are kept and continuously
used. 

It is good and common practice for supervised learning tasks to perform
\emph{model selection and validation} to estimate the out-of-sample predictive
performance and obtain a model that ultimately has good \emph{prediction
accuracy} and generalizes well to unseen data. To perform model selection and
validation, we will use the MSE as a performance metric to evaluate the
performance of a trained model.

To select an appropriate FNN or $k$-NN model and estimate its performance on
unseen data, we will perform \emph{hyperparameter tuning} through an exhaustive
grid search and use $k$-fold CV (note that this $k$ is
different from $k$ in $k$-NN) on a dataset. For a hyperparameter
configuration, $k$-fold CV consists of the following steps: (i)
Divide the dataset into $k$ non-overlapping groups of approximately equal size,
then (ii) save a single fold as a validation dataset and use an ML algorithm to
train a model on the remaining $k - 1$ folds. Lastly, (iii) assess the
performance of the trained model by calculating the MSE on the held-out fold.
Steps (i) and (ii) are repeated $k$ times where a new fold is treated as the
validation dataset. Through this procedure, we obtain error estimates $MSE_1,
MSE_2, \dots, MSE_k$ that can be used to calculate the overall $k$-fold CV score:
\begin{align}
      CV_{k} = \frac{1}{k} \sum_{i = 1}^{k} MSE_i
\end{align}
The $CV_k$ score estimates the error we can expect when using a trained model to
predict new values given unseen data. The model with the best score is chosen in
a model selection context. The value of $k$ is commonly set to 5 or 10, usually
due to computational considerations \citep{jam2013}.

\subsubsection{Optimization}

We can then solve the outer-level optimization problem in
Eq.~\eqref{eqn:1_estimation_problem} with $\hat{\Gamma}$ to model the
relationship in Eq.~\eqref{eqn:sumo_dta_map}. So instead of obtaining arc counts
by performing a time-consuming DTA with SUMO, we get arc counts by querying the
ML model $\hat{\Gamma}$. Since we use nonparametric ML algorithms, the OF is not
explicitly available. Hence, we need to resort to derivative-free methods.
However, now evaluating the OF has become fast, and we can use algorithms that
are less parsimonious in requiring this.

We can thus apply an off-the-shelf global optimization algorithm called the
\emph{basinhopping} algorithm, described in \cite{Dav1997} and implemented by
the SciPy Python library \citep{2020SciPy}. The basinhopping algorithm is
similar to the well-known \emph{simulated annealing} optimization algorithm in
the sense that it consists of two main components that are applied iteratively:
(i) A random perturbation of a current point and (ii) a criterion for accepting
or rejecting the perturbed point based on its OF value.

The basinhopping algorithm extends the idea of the simulated annealing algorithm
by applying a local optimization algorithm to the perturbed point to find a
local minimum of the OF before deciding to accept or reject the point that
resulted in the local minimum. This process is repeated until a maximum number
of iterations has been reached. The local optimization algorithm used in the
basinhopping algorithm as a subprocedure is, in this work, chosen to be a
Sequential Least-Squares Programming (SLSQP) algorithm also implemented by the
SciPy library. This choice is because this algorithm can handle not only the
upper and lower bound constraints in Eq.~\eqref{eqn:2_estimation_problem} but
also the inequality constraints in Eq.~\eqref{eqn:generation_constraints}.

Ultimately, the basinhopping algorithm is applied to find an approximate minimum
of the DODE problem. SUMO will then evaluate the final demand estimate as a
final result.

\section{Experimental Setup} \label{sec:experimental_setup}

We decided to test and compare the methods described in a controlled environment
by synthetically generating an artificial traffic network. To generate this test
data, three primary components are needed: (i) a traffic network where possible
OD pairs $W$ have been identified, (ii) a discretization of the analysis period
$T$, and (iii) a \emph{ground-truth} vector of demands based on (i) and (ii).
The specific parameters used in each of the experiments are described in
Section~\ref{sec:results}. Here, we focus on details of components (i)-(iii).

\begin{figure}[tb]
    \centering
    \subfigure[The network used in the case study. All the nodes on the perimeter of the network can be identified as origins and destinations.]{%
        {
            
            \raisebox{1.25cm}{\includegraphics[trim={2.5cm 2.5cm 2.5cm 2.5cm},clip,keepaspectratio,width=0.38\textwidth]{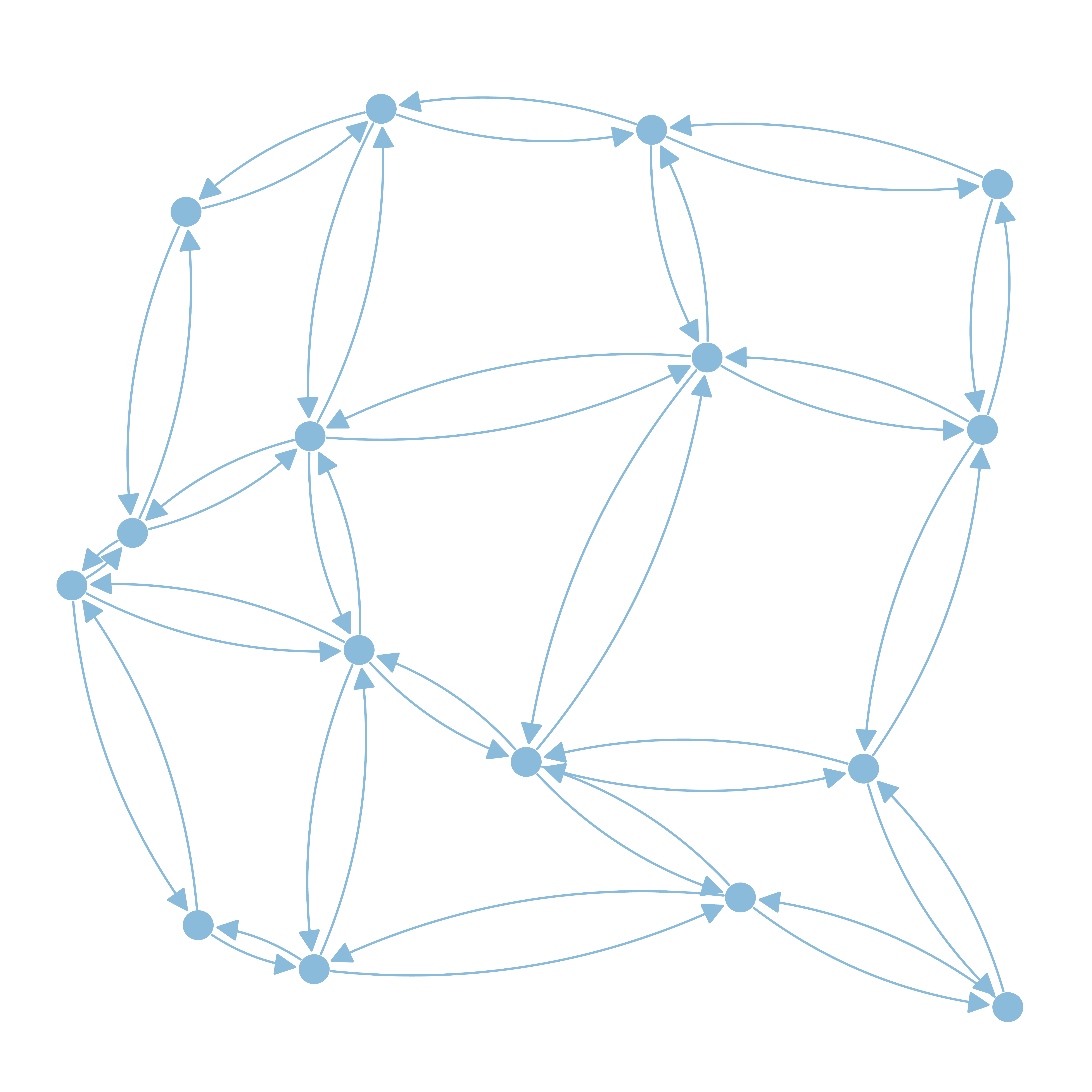}}
            \label{fig:network_illustration}        
        }
    }
    \hspace{1pt}
    \subfigure[Characteristics of the synthetically generated test data using the artificial network.]{%
        {
            \includegraphics[trim={0.5cm 0.5cm 0.25cm 0.5cm},clip,keepaspectratio,width=0.543\textwidth]{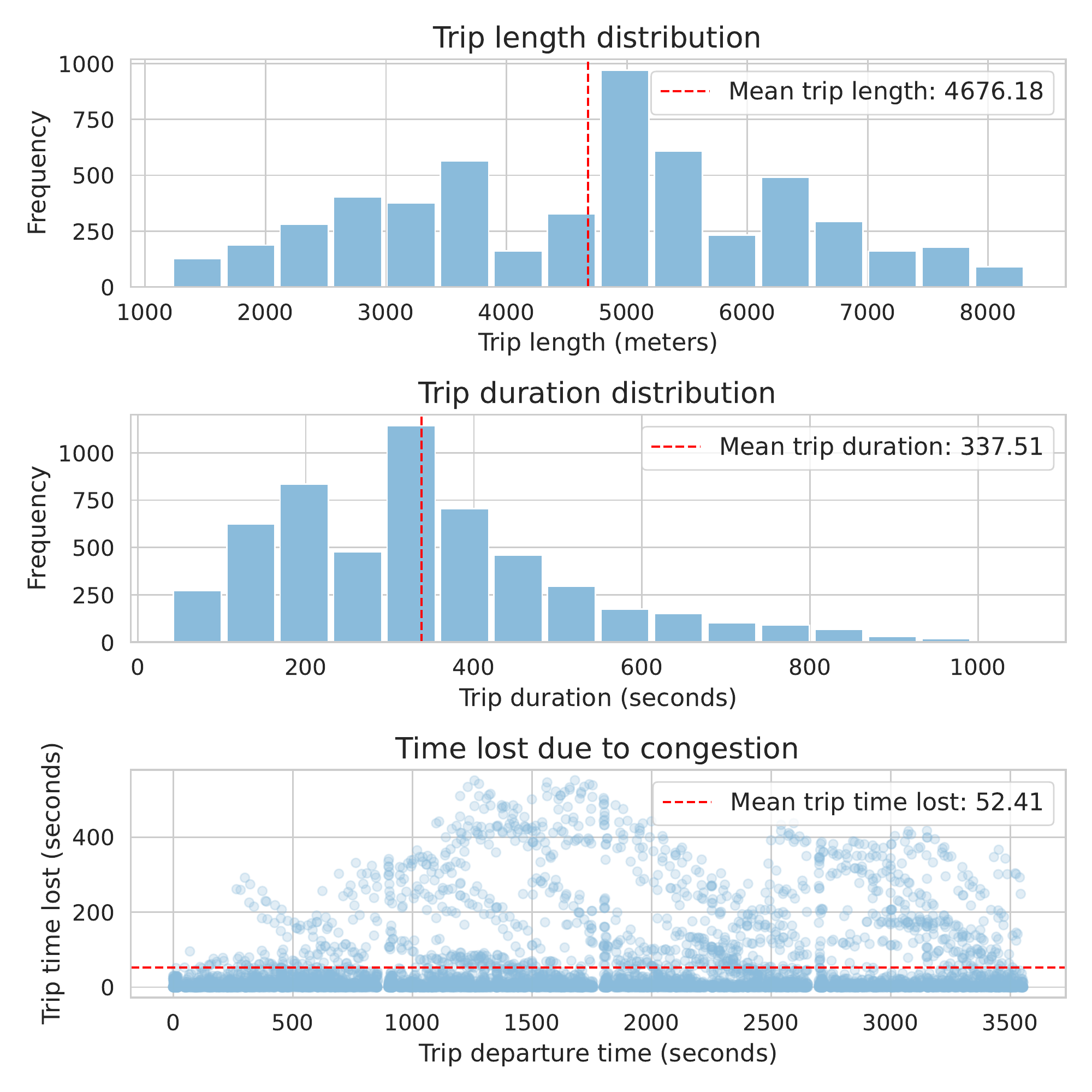}
            \label{fig:network_1_characteristics}        

        }
    }
    \caption{The network used in the case study along with the characteristics of the synthetically generated test data.}
\end{figure}

The network can be defined as an \emph{irregular grid network}. The network was
generated to contain features that can be encountered in a real traffic network.
More specifically, some paths in the network between origins and destinations
are shorter than others in terms of free-flow travel time. The network used is
depicted in Figure~\ref{fig:network_illustration} and described in more detail
in the following.

The network has $48$ directed arcs and $16$ nodes, either representing
intersections or origins and destinations. More precisely, $12$ of the nodes
represent origins and destinations, which means that $O = D$ and $|O| = |D| =
12$, resulting in $132$ OD pairs. Furthermore, the arcs in the network were
prescribed a speed limit of $50 \;
\nicefrac{\textnormal{kilometer}}{\textnormal{hour}}$ and assigned a default
length of $1250$ meters. The coordinates of the nodes in the network are then
perturbed by a uniform random number $\mathcal{U}(0, 1250)$, resulting in arcs
between nodes with different lengths, yielding the irregular grid.

We settled for a one-hour scenario with an analysis period $T = \left[0,
t_{\textnormal{end}}\right)$, where $t_{\textnormal{end}} = 3600 \textnormal{
seconds} = 60 \textnormal{ minutes} = 1 \textnormal{ hour}$. This one-hour time
period was discretized into $4$ time intervals of equal duration, i.e., we have
$S = \left\{1,...,n_S = 4\right\}$, where each of the indices in this set refers
to a time interval of $900 \textnormal{ seconds} = 15 \textnormal{ minutes}$
duration. Given this information and a network with OD pairs $W$, a ground-truth
vector of demands is constructed by sampling the number of trips defined for
each variable according to a continuous uniform distribution with lower bound
$a$ and upper bound $b$:
\begin{align}
    x_{ws}^{\textnormal{True}} \sim U(a, b), \; \forall w \in W, s \in S.
    \label{eqn:ground_truth}
\end{align}
From the ground-truth vector of demands, trip productions $o_i$ for $i =
1, \dots, |O|$ that enter into the generation constraints (in
Eq.~\eqref{eqn:generation_constraints}) were obtained in addition to upper and
lower bounds on the estimated demands:
\begin{align}
    \mathbf{x}_{l} = \mathbf{0}  \;\; \textnormal{and} \;\;  \mathbf{x}_{u} = 1.5 \cdot \max\left(\mathbf{x}\right) \cdot \mathbf{1},
\end{align}
where $\mathbf{0}$ is an $m \cdot n_S$ vector of zeros and $\mathbf{1}$ is an $m
\cdot n_S$ vector of ones.

Finally, the ground-truth vector of demands was given to SUMO as input
and the iterative DTA procedure implemented by SUMO was run for $15$ iterations
yielding the vector $\hat{\mathbf{c}}$ of observed arc counts.

For the seed demand vector $\tilde{\mathbf{x}} \in \mathbb{R}^{m \cdot
n_S}$ we consider two scenarios:
\begin{description}
    \item[Prior low-demand vector (LD):] This scenario simulates the case where 
    lower demand estimates are available from a previous study:
    \begin{align}
        \tilde{x}_{ws} =  x_{ws}^{\textnormal{True}} \cdot \left(0.7 + 0.3 \cdot u_{ws}\right), \qquad u_{ws} \sim U(0, 1), \; \forall w \in W, \forall s \in S.
    \end{align}{}
    \item[Prior high-demand vector (HD):] This scenario simulates the case where 
    higher demands are available from a previous study:
    \begin{align}
        \tilde{x}_{ws} =  x_{ws}^{\textnormal{True}} \cdot \left(0.9 + 0.3 \cdot u_{ws}\right), \qquad u_{ws} \sim U(0, 1), \; \forall w \in W, \forall s \in S.
    \end{align}{}
\end{description}

These two scenarios were constructed based on the guidelines provided in
\cite{ant2016} and the computational experiments performed in \cite{cip2011,
mas2018}.

The ground-truth vector of demands was generated according to
Eq.~\eqref{eqn:ground_truth}, using lower bound $a = 1$ and upper bound $b =
20$. The value for the upper bound $b$ depends on the number of OD pairs and the
size of the network, i.e., the network supply, which essentially determines how
much traffic the network can handle at a given point in time. In this case, the
appropriate upper bound was found experimentally by sampling different
ground-truth vectors of demand for increasing values of $b$. A DTA was
performed by SUMO for each of the generated ground-truth vectors, and the
resulting output from the DTA was assessed by checking the level of congestion
present in the network. Finally, the ground-truth vector of demands with
corresponding upper bound $b$ that generated a scenario with light congestion was
chosen as the basis for the case study using this network.

The scenario that was deemed the most appropriate was determined based on
assessing the characteristics of the three plots in
Figure~\ref{fig:network_1_characteristics}. In the first plot, we observe that
the average trip length is $4.6$ kilometers and, in the next plot, that the
average travel time for a trip is $338$ seconds. Finally, keeping these two
attributes in mind, we observe in the last plot that an average delay of $52$
seconds ($\nicefrac{52 \textnormal{ seconds}}{338 \textnormal{ seconds}} =
0.15\%$) of the total travel time is due to congestion and driving below the
possible free-flow speed. Furthermore, some unfortunate vehicles experience a
delay far more significant than the average of $52$ seconds, while on the other
hand, most vehicles experience only little to no travel time delays. Therefore,
based on these scenario attributes, we classify the scenario as a scenario with
light congestion.

For the results we present in the following sections, it is assumed that sensors
cover $100$\% of the network. A good overview of the problem of determining
appropriate locations to place sensors in a traffic network, along with some
recent developments on this problem, is given in \cite{vit2014, hu2014,
sal2019}. In addition to this, the arc counts obtained by the sensors are
assumed to be error-free. This is, of course, not the case in a real traffic
network, but it simplifies the analysis and comparison of the different solution
approaches.  Lastly, we mention that all the experiments were executed using an
initial feasible solution consisting of all ones, i.e., $\mathbf{x}_{\tau = 0} =
\mathbf{1}$. This initial guess was chosen so to not confine the search for a
minimum to the vicinity of the seed demands.

To compare the performance of the different approaches, we report the Root Mean
Squared Error (RMSE). If we consider place-holders $\mathbf{y}, \hat{\mathbf{y}}
\in \mathbb{R}^{n}$, where $\mathbf{y}$ is a vector of observed values and
$\hat{\mathbf{y}}$ is a vector of estimated values, then the RMSE value is
defined as follows:
\begin{align}
   \textnormal{RMSE}(\mathbf{y}, \hat{\mathbf{y}}) = \left(\frac{1}{n}\sum_{j = 1}^{n} (y_{i} - \hat{y}_{i})^2\right)^{\frac{1}{2}}.
\end{align}

This error measure is easily interpretable because the error is given in the
same units as the data. Furthermore, the RMSE is reported separately for each of
the different quantities entering the optimization problem, i.e., we report the
RMSE statistic for the arc counts and the demands. Finally, we remark
that for the demands, the RMSE is calculated based on the ground-truth
values and not on the seed demand vectors that enter the optimization
problem.

We note that the size of the network and the number of OD pairs, arcs
equipped with sensors and time intervals of the dataset described are
similar to those used in other algorithm comparison studies
\citep{cip2011, car2017, nig2018_1}, while larger sizes and numbers are
used in \cite{tol2013, ant2016, dju2017, mas2018}.

\section{Experimental Results} \label{sec:results}

We implemented all methods in Python with external calls to the SUMO simulator.
The primary Python libraries used include Pandas \citep{mck2010} for data
organization purposes, NumPy \citep{harris2020array} for its numerical methods
for fitting polynomials, along with SciPy \citep{2020SciPy} for its off-the-shelf
global optimization algorithms and scikit-learn's \citep{sklearn} implementation
of FNN and $k$-NN used in the ML approach. A snapshot
of the repository \citep{github_demand_estimation_data_available} containing the
codebase and the data supporting this study's findings are archived and openly
available in Zenodo \citep{zenodo_demand_estimation_data_available}. All
experiments described in this section were performed on a machine with an AMD
Ryzen 7 1800X eight-core processor and 32 GB of RAM, running Manjaro Linux.

% \paragraph{Computational Budget}
\subsection{Computational Budget}
\begin{figure}[tb]
  \centering
  \includegraphics[width = 0.75\textwidth]{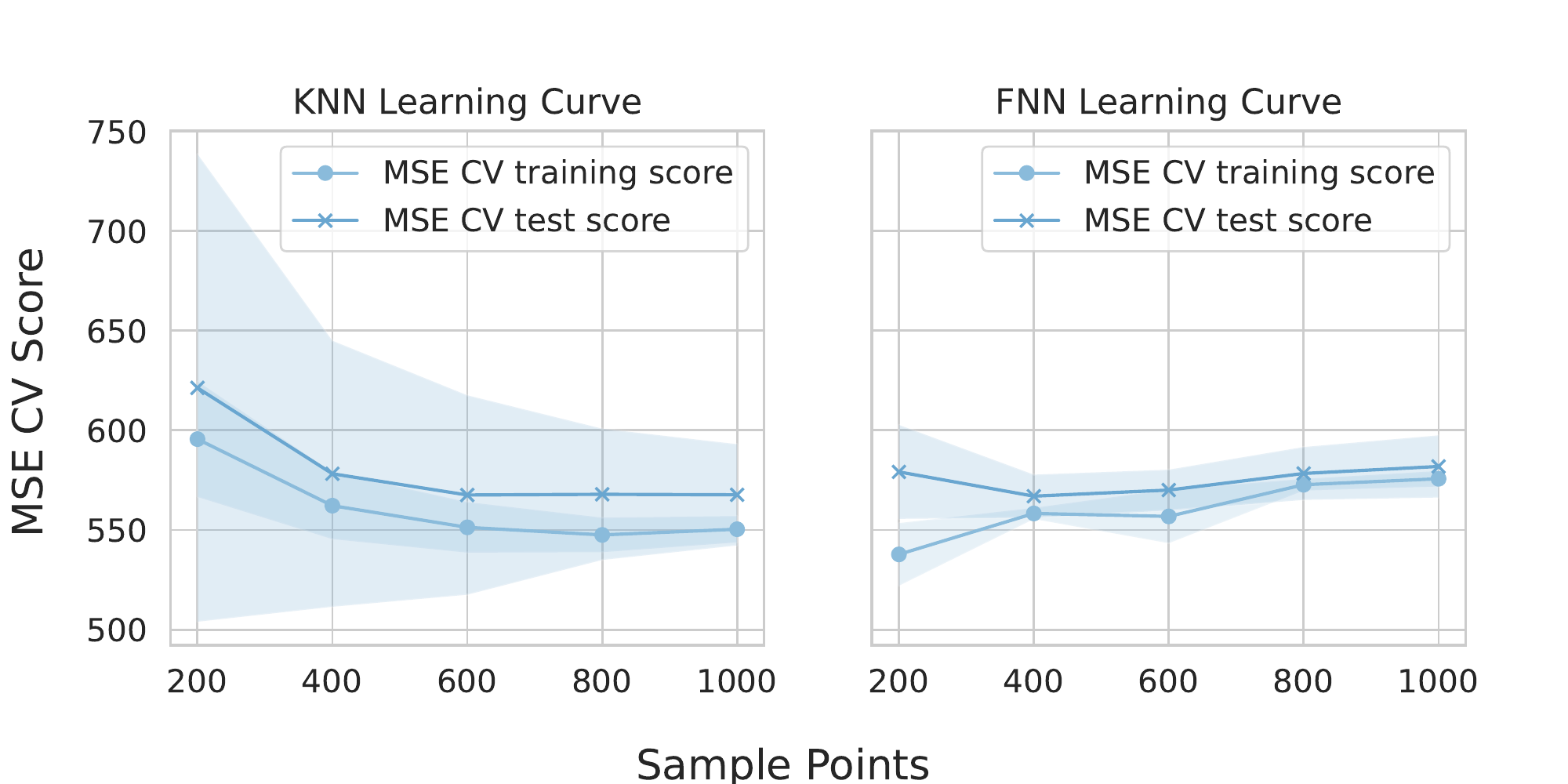}
  \caption{
    The learning curves and CV results associated with the training of the $k$-NN and FNN model.
    The curves show the mean CV and training scores, while the shaded bands show the standard deviation.
    The datapoints in the curves are related to the best mean CV score obtained for different sample sizes
    and over all tested hyperparameter configurations.    
    }
  \label{fig:learning_curve}
\end{figure}

We studied the ability of the FNN and $k$-NN ML algorithms to learn the
input-output relationship using a different number of sample points and
hyperparameter configurations while the CV scores were monitored. The resulting
\emph{learning curves} are displayed in Figure~\ref{fig:learning_curve}. We
observe that by increasing the size of the training set, the CV scores improve
but only slightly.

As we want a model that performs as well as possible on unseen data (the test
scores) and is parsimonious in the usage of SUMO to perform DTAs, we settled to
use $200$ points for the dataset to train the ML approaches. This means
that we run SUMO for about $200$ DTAs, which took $3912.92 \textnormal{ sec} =
65.21\textnormal{ min}$ on our computational environment. We note that an
additional DTA (the $201$st DTA) will be needed by the ML approaches to evaluate
a final estimate. For the sake of fairness, we thus set a soft upper bound on
the number of allowed OF evaluations to $201$ for all other
approaches\footnote{Due to the inherent nature of some of the applied
approaches, a few more OF evaluations and thus DTAs might have
occurred in an iteration and the total number of OF evaluations
might not sum up to $201$ exactly.}.

% \paragraph{Parameter Tuning}
\subsection{Parameter Tuning}

The gradient-based approaches, i.e., the AM-based and AM-free approach, do not
need parameter tuning before tackling the DODE problem. Indeed, the only
parameter, the step length, is determined automatically by the solution of the
line search in Eq.~\eqref{eqn:line_search}. Thus, these algorithms are
straightforward to apply, and no specific consideration has to be made
concerning their parameters.

For the ML approaches, the FNN and the $k$-NN model, the hyperparameters were
fine-tuned using a grid search coupled with $5$-fold CV using the $200$
generated sample points that act as re-usable historical simulation data.

The grid search was performed on the values reported in
Table~\ref{tab:hyperparameters}. The best performance was achieved by the
following settings:
\begin{itemize}
    \item The FNN model with hidden layers containing
    units $(0.75\cdot m \cdot n_S, 0.50\phantom{0} \cdot m \cdot n_S,
    0.50\phantom{0} \cdot m \cdot n_S) = (l_1=396, l_2=264, l_3=264)$ along with
    regularization parameter $\lambda=0.01$. This trained model is used in
    method M3 in Table~\ref{tab:summary_table}. This model achieved a $5$-fold
    CV score of $CV_5 = 579.1282$.

 \item The $k$-NN model using $k = 43$ neighbours. This trained model is used in
    method M4 in Table~\ref{tab:summary_table}. This model achieved a $5$-fold
    CV score of $CV_5 = 621.2651$.
\end{itemize}

For the subsequent optimization step using a trained ML model (FNN or $k$-NN) as
a surrogate model for DTA, the basinhopping algorithm was used with default
parameters. Also, the local SLSQP minimizer was used with default parameters.
Preliminary experiments showed that $10$ iterations of the basinhopping
algorithm were enough to reach convergence and that improvements after that
stage become rare.
% SUMO, in the end, evaluates the estimate obtained by the
%optimization step through a DTA as the final $201$st OF evaluation.

\subsection{Quality of Estimates}

In Table~\ref{tab:summary_table}, we present an overview of all
experiments conducted together with summary results, that we detail in the following.

\begin{table}
    \footnotesize
    \centering
      \caption{ 
          An overview and summary of the results obtained in the different
          experiments E1-E12. Here, we define M1: AM-based approach,
          M2: AM-free approach, M3: FNN predicting $\Gamma$,
          M4: $k$-NN predicting $\Gamma$. For methods M3-M4 only the
          running time of the optimization step is reported. 
      }
      \setlength{\tabcolsep}{2pt}
\begin{tabular}{c c c c c c c c}
    \toprule
    \multicolumn{1}{c}{Method} & \multicolumn{1}{c}{Experiment} & \multicolumn{1}{c}{Seed matrix} &
    % \multicolumn{1}{c}{Total OF Evaluations} & \multicolumn{1}{c}{Running Time (sec)} &
    \multicolumn{1}{c}{OF evaluations} & \multicolumn{1}{c}{Running time (sec)} &
    \multicolumn{1}{c}{RMSE$(\mathbf{x}, \mathbf{x}^{\textnormal{True}})$} &
    \multicolumn{1}{c}{RMSE$(\mathbf{c}(\mathbf{x}), \hat{\mathbf{c}})$} \\

    \cmidrule(lr){1-1} \cmidrule(lr){2-2}
    \cmidrule(lr){3-3} \cmidrule(lr){4-4}
    \cmidrule(lr){5-5} \cmidrule(lr){6-6}
    \cmidrule(lr){7-7} \cmidrule(lr){8-8}
    % METHOD              & EXPERIMENT          % SEED MATRIX             % OF EVALUATIONS     % Running time     % RMSE ODMAT        % RMSE COUNTS
    \multirow{3}{*}{M1} & \multirow{1}{*}{E1} & \multirow{1}{*}{None} & \multirow{1}{*}{201} & \multirow{1}{*}{2836.4469} & \multirow{1}{*}{ 1.3413 } & \multirow{1}{*}{9.4626} \\
                        & \multirow{1}{*}{E2} & \multirow{1}{*}{LD}   & \multirow{1}{*}{201} & \multirow{1}{*}{2845.5479} & \multirow{1}{*}{ 1.3413 } & \multirow{1}{*}{9.4626} \\
                        % & \multirow{1}{*}{E3} & \multirow{1}{1cm}{MD}   & \multirow{1}{*}{27} & \multirow{1}{*}{ 3.07} & \multirow{1}{*}{14.50} \\
                        & \multirow{1}{*}{E3} & \multirow{1}{*}{HD}   & \multirow{1}{*}{201} & \multirow{1}{*}{2770.9962} & \multirow{1}{*}{ 1.2516 } & \multirow{1}{*}{9.6358} \\
    \midrule

    % METHOD              & EXPERIMENT          % SEED MATRIX             % OF EVALUATIONS    % Running time      % RMSE ODMAT        % RMSE COUNTS
    \multirow{3}{*}{M2} & \multirow{1}{*}{E4} & \multirow{1}{*}{None} & \multirow{1}{*}{203} & \multirow{1}{*}{2552.5320} & \multirow{1}{*}{9.0188} & \multirow{1}{*}{19.4008} \\
                        & \multirow{1}{*}{E5} & \multirow{1}{*}{LD}   & \multirow{1}{*}{201} & \multirow{1}{*}{2366.2931} & \multirow{1}{*}{ 2.0342} & \multirow{1}{*}{13.5121} \\
                        % & \multirow{1}{*}{E7} & \multirow{1}{1cm}{MD}   & \multirow{1}{*}{6375} & \multirow{1}{*}{ \textbf{2.58}} & \multirow{1}{*}{13.89} \\
                        & \multirow{1}{*}{E6} & \multirow{1}{*}{HD}   & \multirow{1}{*}{202} & \multirow{1}{*}{2475.8943} & \multirow{1}{*}{ 6.7992} & \multirow{1}{*}{16.1219} \\
    \midrule
    % METHOD              & EXPERIMENT          % SEED MATRIX             % OF EVALUATIONS    % Running time      % RMSE ODMAT        % RMSE COUNTS
    \multirow{3}{*}{M3} & \multirow{1}{*}{E7}& \multirow{1}{*}{None} & \multirow{1}{*}{201 } & \multirow{1}{*}{110.9053} & \multirow{1}{*}{9.6494} & \multirow{1}{*}{20.8409} \\
                        & \multirow{1}{*}{E8}& \multirow{1}{*}{LD}   & \multirow{1}{*}{201 } & \multirow{1}{*}{116.7221} & \multirow{1}{*}{9.6338} & \multirow{1}{*}{20.4193} \\
                        % & \multirow{1}{*}{E19}& \multirow{1}{1cm}{MD}   & \multirow{1}{*}{10004 + 11} & \multirow{1}{*}{12.64} & \multirow{1}{*}{27.45} \\
                        & \multirow{1}{*}{E9}& \multirow{1}{*}{HD}   & \multirow{1}{*}{201} & \multirow{1}{*}{107.7270} & \multirow{1}{*}{9.6493} & \multirow{1}{*}{20.8409} \\
    \midrule
    % METHOD              & EXPERIMENT          % SEED MATRIX             % OF EVALUATIONS      % Running time           % RMSE ODMAT        % RMSE COUNTS
    \multirow{3}{*}{M4} & \multirow{1}{*}{E10} & \multirow{1}{*}{None} & \multirow{1}{*}{201 } & \multirow{1}{*}{103.1716} & \multirow{1}{*}{9.6728} & \multirow{1}{*}{21.1067} \\
                        & \multirow{1}{*}{E11}& \multirow{1}{*}{LD}   & \multirow{1}{*}{201} & \multirow{1}{*}{109.7821} & \multirow{1}{*}{9.6338} & \multirow{1}{*}{20.4193} \\
                        % & \multirow{1}{*}{E11}& \multirow{1}{1cm}{MD}   & \multirow{1}{*}{10004 + 7 } & \multirow{1}{*}{12.69} & \multirow{1}{*}{27.99} \\
                        & \multirow{1}{*}{E12}& \multirow{1}{*}{HD}   & \multirow{1}{*}{201} & \multirow{1}{*}{102.5330} & \multirow{1}{*}{9.6339} & \multirow{1}{*}{20.4193} \\

    \bottomrule
\end{tabular}

      \label{tab:summary_table}    
\end{table}

% \subparagraph{Gradient-based Approaches}
\subsubsection{Gradient-based Approaches}

\begin{figure}[tb]
    \centering
        \subfigure[The ground-truth values plotted against the estimated values.]{
            \includegraphics[trim={0.225cm 0.1cm 1.25cm 0.0cm},clip,width=0.59\linewidth,keepaspectratio]{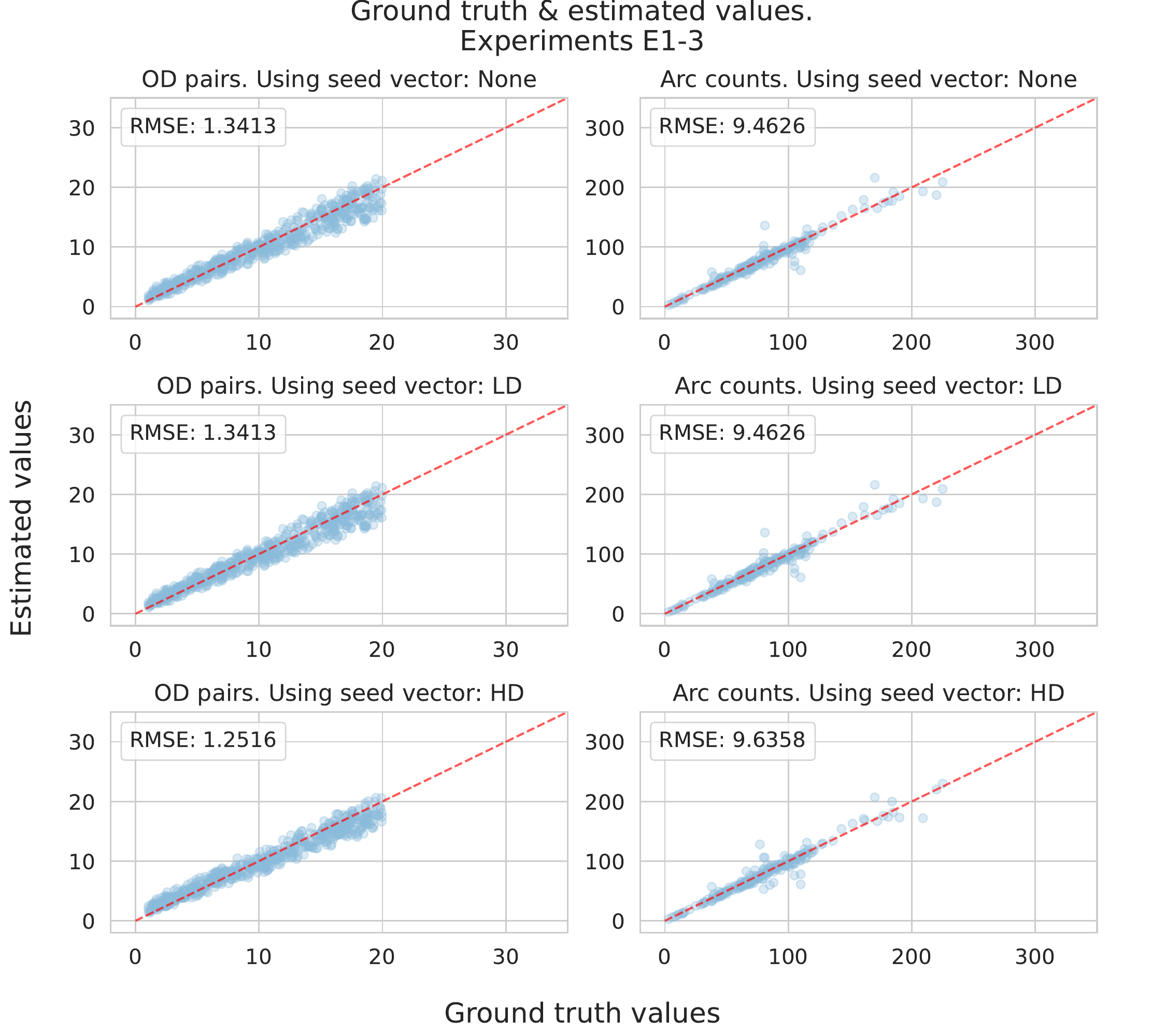}
            \label{fig:assignmat_approach_estobs}
        }
        \hspace{1.0pt}
        \subfigure[The OF value history. Utilizing output data directly from SUMO $\Gamma(\mathbf{x})$.]{
            \includegraphics[trim={0.0cm 0.0cm 0.0cm 0.0cm},clip,width=0.35\linewidth]{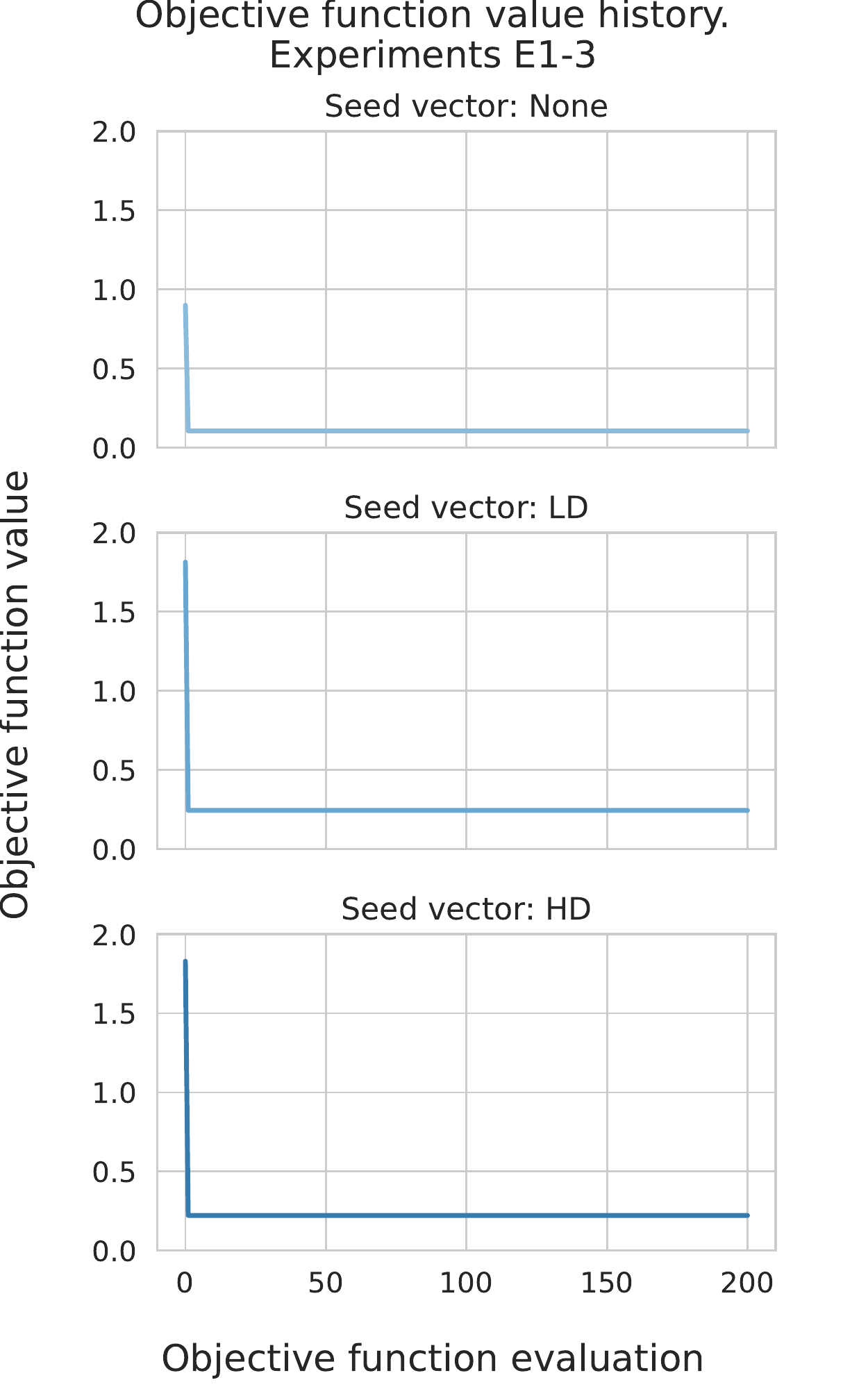}
            \label{fig:assignmat_approach_iters}
        }
        \caption{The final results from experiments E1-E3, that used the gradient and AM-based approach.}
        \label{fig:assignmat_approach}
\end{figure}

\begin{figure}[tb]
    \centering
    \subfigure[The ground-truth values plotted against the estimated values.]{
        \includegraphics[trim={0.225cm 0.1cm 1.25cm 0.0cm},clip,width=0.59\linewidth,keepaspectratio]{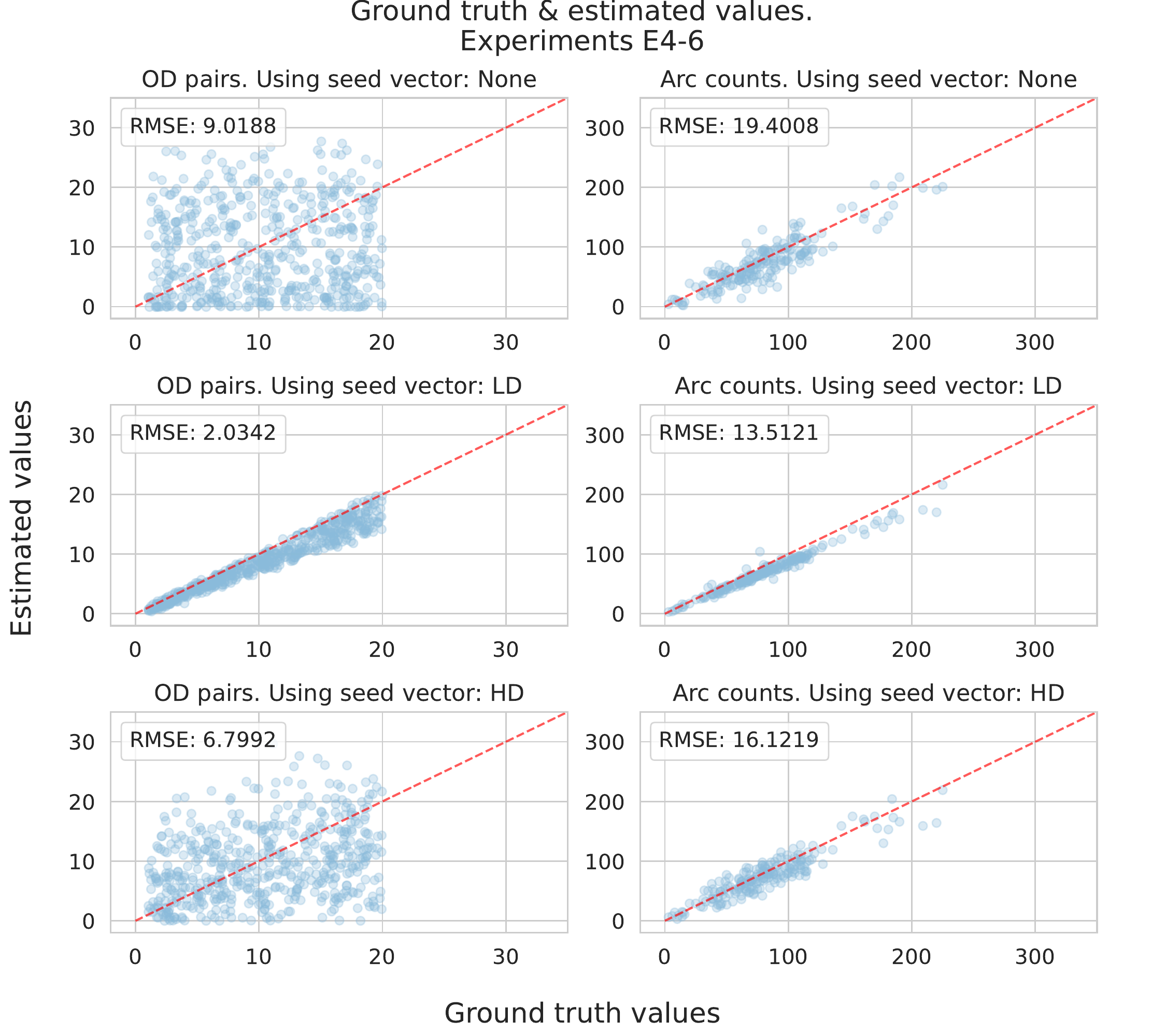}
        \label{fig:spsa_approach_estobs}
    }
    \hspace{1.0pt}
    \subfigure[The OF value history. Utilizing output data directly from SUMO $\Gamma(\mathbf{x})$.]{
        \includegraphics[trim={0.0cm 0.0cm 0.0cm 0.0cm},clip,width=0.35\linewidth]{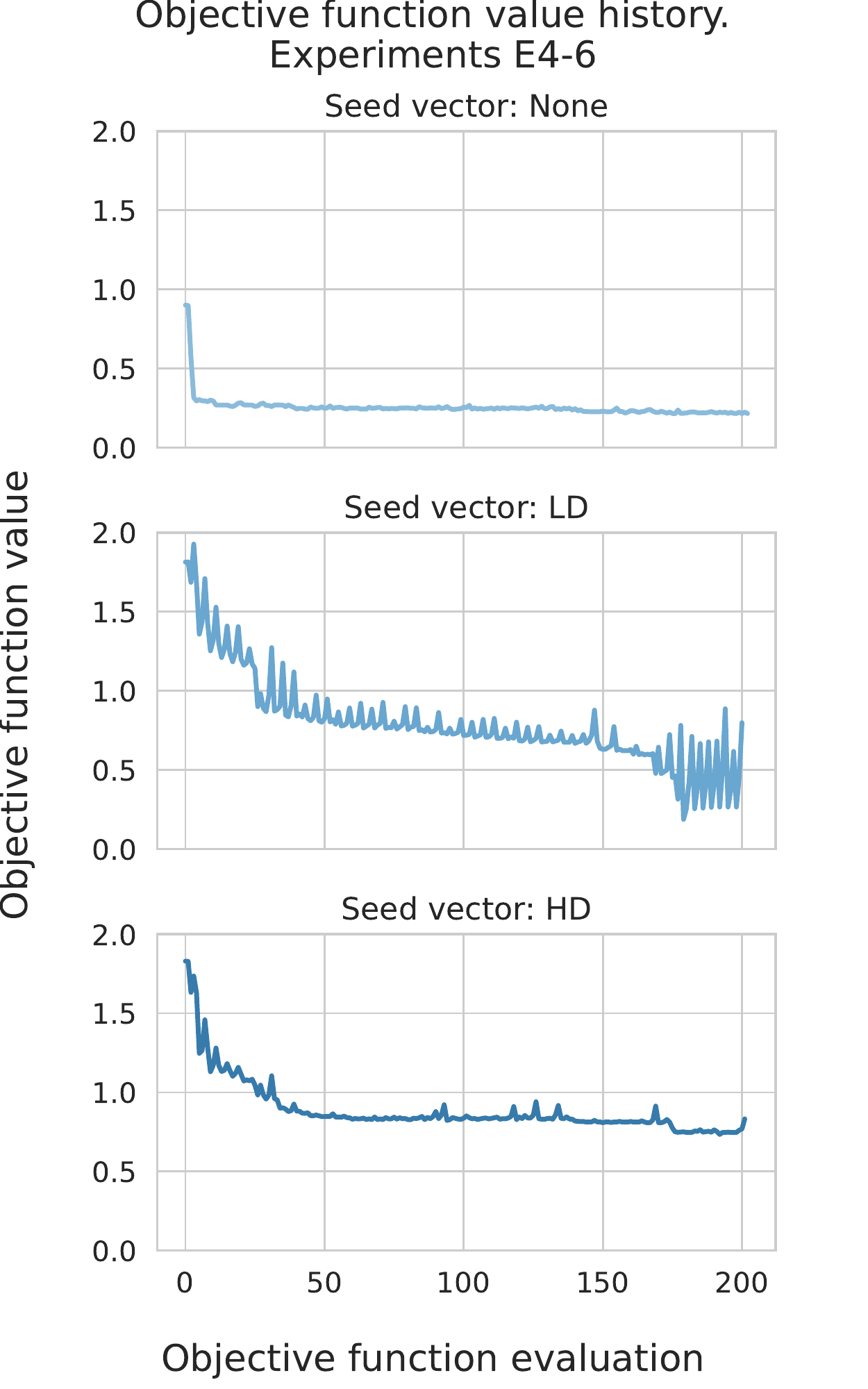}
        \label{fig:spsa_approach_iters}
    }
    \caption{The final results from experiments E4-E6, that used the gradient and AM-free approach.}
    \label{fig:spsa_approach}
\end{figure}

The final results obtained from the gradient-based approaches are plotted in
Figure~\ref{fig:assignmat_approach_estobs} and
Figure~\ref{fig:spsa_approach_estobs}. In these figures, the estimated
quantities $\mathbf{x}$ and $\mathbf{c}(\mathbf{x})$ have been plotted against
the ground-truth vector of demands $\mathbf{x}^{\textnormal{True}}$ and
the corresponding vector of observed arc counts $\mathbf{c}(\mathbf{x})$,
respectively. The corresponding plots that show the OF value history of the two
algorithms are shown in Figure~\ref{fig:assignmat_approach_iters} and
Figure~\ref{fig:spsa_approach_iters}.

If we compare the results in Figure~\ref{fig:assignmat_approach_estobs} and
Figure~\ref{fig:spsa_approach_estobs} between the two gradient-based approaches,
we observe that the AM-based approach outperforms the AM-free approach. In fact,
the AM-based approach is able to consistently overcome local minima and reach a
good demand estimate. This can be seen in
Figure~\ref{fig:assignmat_approach_estobs} where the estimated demands and the
ground-truth vector of demands align well, and so do also the estimated arc
counts and the observed arc counts.

If we look at the progress of the two algorithms in
Figure~\ref{fig:assignmat_approach_iters} and
Figure~\ref{fig:spsa_approach_iters}, we see that the AM-based
approach can find a good estimate within only a single OF
evaluation, while the AM-free approach uses a much larger number
of evaluations to arrive at good estimates.

If we return to Table~\ref{tab:summary_table} and the results relating to
methods M1 and M2's computational time, we observe that the AM-free (M2) appears
to be faster than the AM-based (M1) approach. This is the case, even though two
OF evaluations and thus DTA computations are needed in each iteration of M1,
while M2 needs up to four. The prolonged computational time can be attributed to
slow gradient calculations due to AM data organization, which takes time in
Python. As the network and scenario are scaled, we expect the case to be
inverted, i.e., data organization will take less time than performing DTA's.


\subsubsection{ML Approaches}

\begin{table}[b]
    \footnotesize
    \centering
    \caption{The ML model hyperparameters chosen for
    hyperparameter tuning. Preliminary experiments were performed to narrow
    down the domains of the hyperparameters. }
    \begin{tabular}{c c c}
        \toprule
        \multicolumn{1}{c}{Machine learning alg.} & \multicolumn{1}{c}{Hyperparameter name} & \multicolumn{1}{c}{Domain} \\
        \cmidrule(lr){1-1} \cmidrule(lr){2-2} \cmidrule(lr){3-3}
        \multirow{4}{*}{FNN} &
        \multirow{3}{*}{Architecture}
            & \multirow{1}{*}{ $\{(0.75\cdot m \cdot n_S, 0.50\phantom{0} \cdot
                m \cdot n_S, 0.50\phantom{0} \cdot m \cdot n_S),$ } \\
            & & \multirow{1}{*}{ $\phantom{\{}(0.75\cdot m \cdot n_S,
                0.25\phantom{0} \cdot m \cdot n_S, 0.25\phantom{0} \cdot m \cdot
                n_S),$ } \\
            & & \multirow{1}{*}{ $\phantom{\{}(0.75\cdot m \cdot n_S, 0.125
                \cdot m \cdot n_S, 0.125 \cdot m \cdot n_S)\} $ } \\
        & \multirow{1}{*}{Regularization $\lambda$}  &
        \multirow{1}{*}{$\{0.0001, 0.001, 0.01, 0.1\}$} \\
        \midrule
        \multirow{1}{*}{$k$-NN} & \multirow{1}{*}{Neighbours $k$} &
        \multirow{1}{*}{$\{\lfloor 2^z \rfloor: z = 2 + i \cdot
        \nicefrac{(6-2)}{14}, \forall i = 0, \dots, 14\}$} \\
        \bottomrule
    \end{tabular}
    \label{tab:hyperparameters}
\end{table}

\begin{figure}[tb]
    \centering
    \subfigure[The ground-truth values plotted against the estimated values.]{
        \includegraphics[trim={0.225cm 0.1cm 1.25cm 0.0cm},clip,width=0.59\linewidth,keepaspectratio]{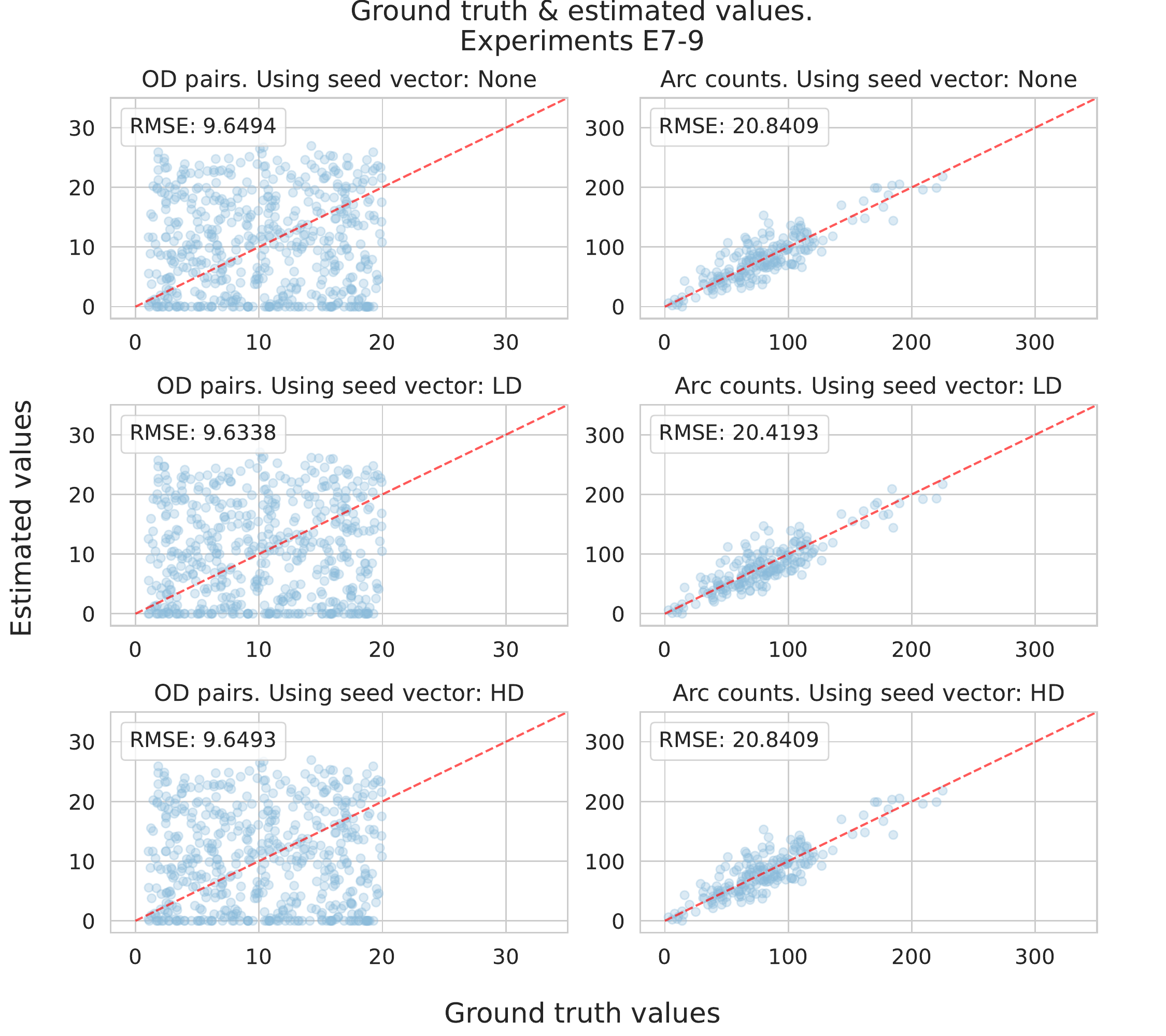}
        \label{fig:knn_approach_estobs}
    }
    \hspace{1.0pt}
    \subfigure[The OF value history. Utilizing output data from the surrogate $\hat{\Gamma}(\mathbf{x})$.]{
        \includegraphics[trim={0.0cm 0.0cm 0.0cm 0.0cm},clip,width=0.35\linewidth]{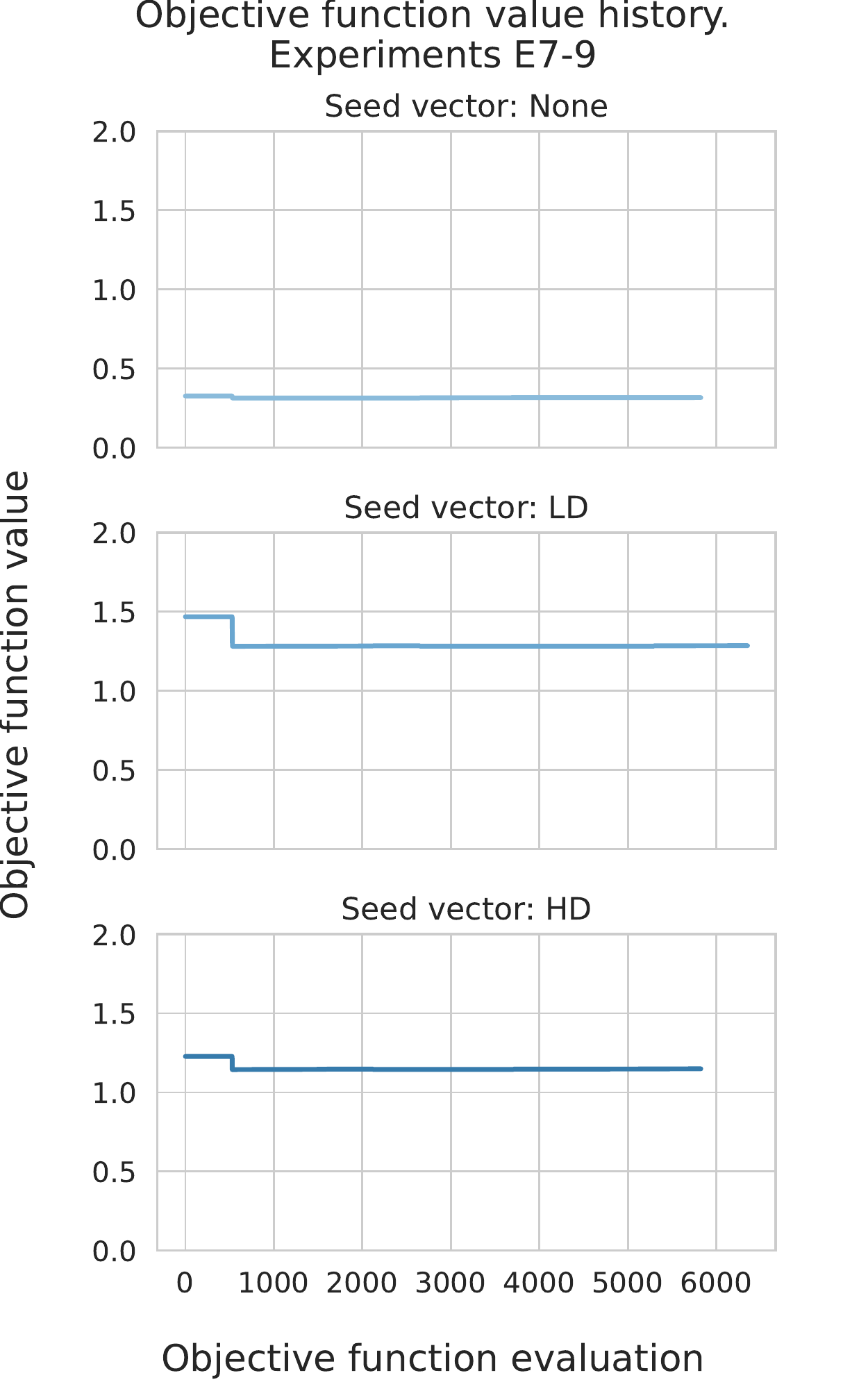}
        \label{fig:knn_approach_iters}
    }
    \caption{The final results from experiments E7-E9, where the $k$-NN model predict $\Gamma$.}
    \label{fig:knn_approach}
\end{figure}

\begin{figure}[tb]
    \centering
    \subfigure[The ground-truth values plotted against the estimated values.]{
        \includegraphics[trim={0.225cm 0.1cm 1.25cm 0.0cm},clip,width=0.59\linewidth,keepaspectratio]{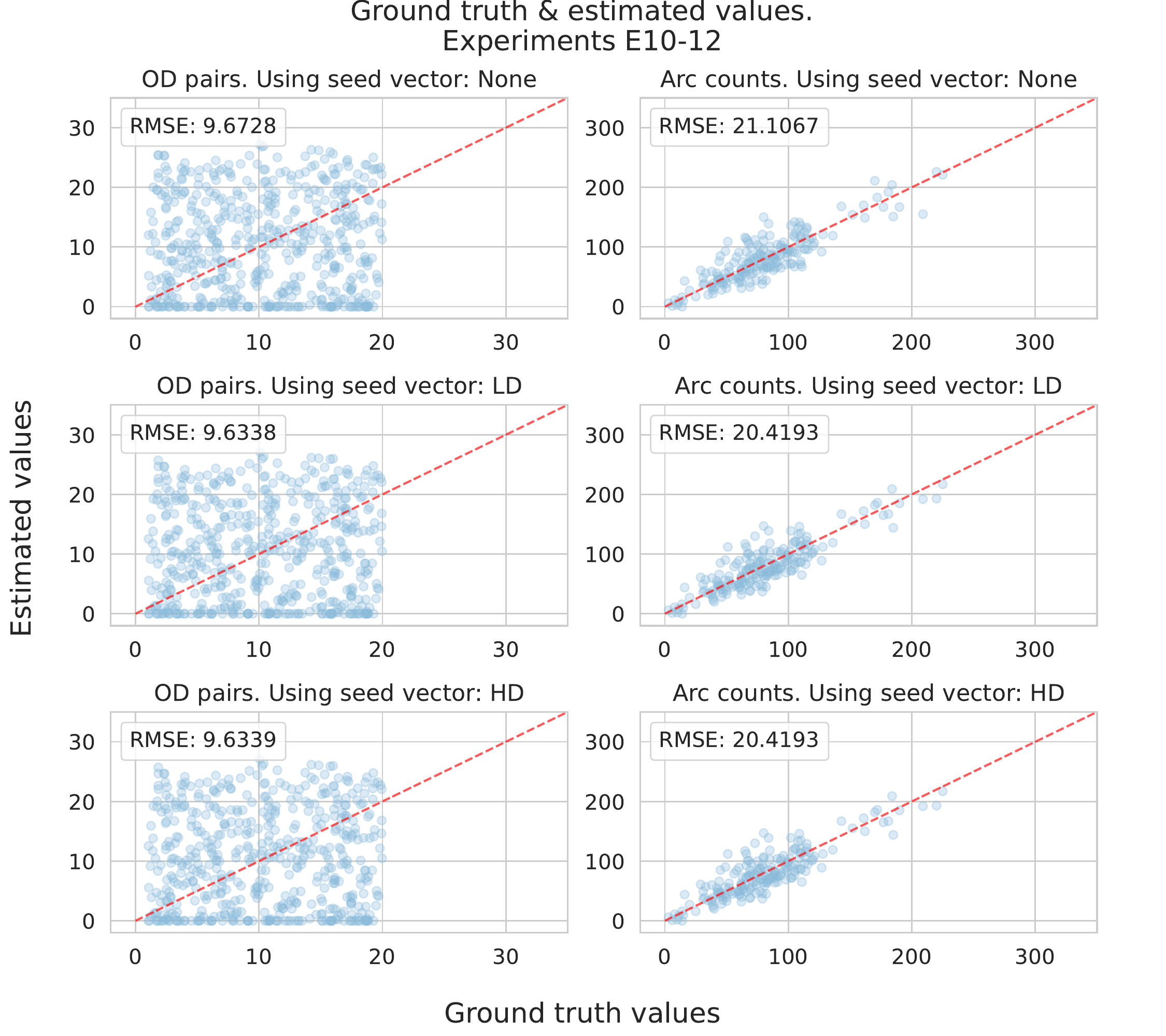}
        \label{fig:fnn_approach_estobs}
    }
    \hspace{1.0pt}
    \subfigure[The OF value history. Utilizing output data from the surrogate $\hat{\Gamma}(\mathbf{x})$.]{
        \includegraphics[trim={0.0cm 0.0cm 0.0cm 0.0cm},clip,width=0.35\linewidth]{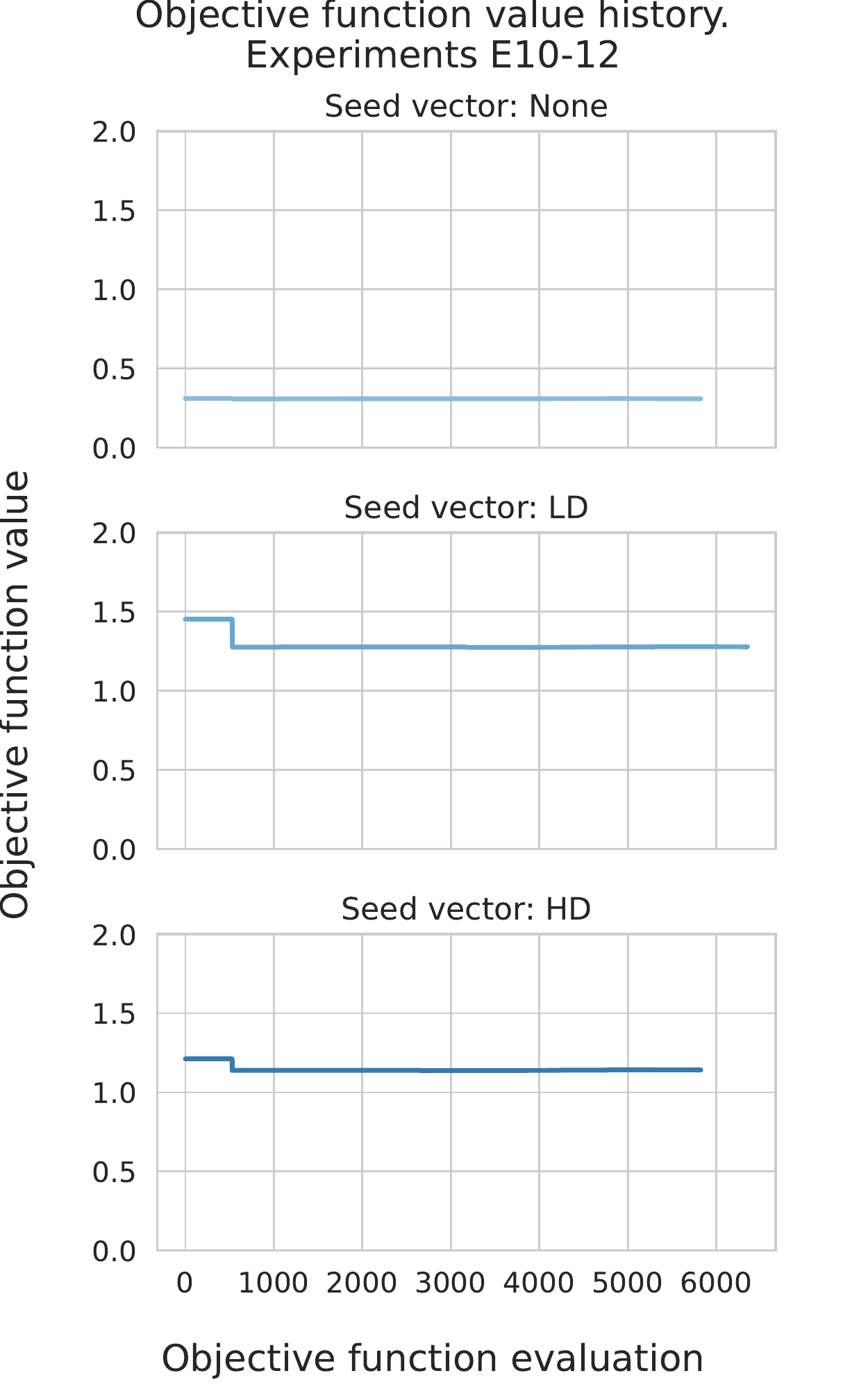}
        \label{fig:fnn_approach_iters}
    }
    \caption{The final results from experiments E10-E12, where the FNN model predict $\Gamma$.}
    \label{fig:fnn_approach}
\end{figure}

The final estimate reached by the $k$-NN approach is shown in
Figure~\ref{fig:knn_approach} and the one reached by the FNN approach in
Figure~\ref{fig:fnn_approach}.  The figures also report the development of
solution quality over the $10$ iterations of the basinhopping algorithm which
typically yield around $6000$ OF evaluations using the surrogate
model $\hat{\Gamma}$.

Across the different experiments with the ML approaches, we observe that the two
ML models achieve very similar results. Under both models, the estimated and
observed arc counts align well but the ground-truth demands and the estimated
demands do not. In other terms, these approaches can find demand estimates that
reproduce the arc counts well but are in no way similar to the ground-truth
demands.

In the optimization step of the ML approaches, we observe a general lack of
progress.  In fact, an examination of the quality of the data points against the
final point after the optimization step show that the ML approaches seem able to
improve only slightly the estimates already available initially.

Furthermore, if the OF values pertaining to the first and the second term of the
OF is inspected, then it appears that the first term (in
Eq.~\eqref{eqn:objective_function_terms}) relating to the demands is not worked
on as much as the second term that takes into account the arc counts. This
behavior helps explain why the ground-truth demands and the estimated demands do
not align well even when a seed vector is provided to guide the search. Another
reason for the overall poor performance can also be attributed to the
goodness-of-fit of the trained ML models. If the trained models do not
sufficiently capture the complex input-output relationship of a DTA performed by
SUMO, then we might only be able to explore or reach local minima of the DODE
problem in the subsequent optimization step that are not actual minima. A hint
that this might be the case is the high CV scores reported in
Figure~\ref{fig:learning_curve}.

% \clearpage

% \paragraph{Scaling}
\subsection{Scaling}

Under the assumption that the training set of data is available initially in the
form of historical data, then performing the optimization step on the ML models
took on average $108.47\textnormal{ sec}$ (average running time of M3-M4 in
Table~\ref{tab:summary_table}). In contrast, the average time for the
gradient-based approaches using SUMO directly for DTA was $2641.28\textnormal{
sec}$ (average running time of M1-M2 in Table~\ref{tab:summary_table}).

Further, depending on the network size, the ML approaches could scale better
than the gradient-based approaches. Compared to the gradient-based approaches
that perform many DTAs, starting from scratch in each new optimization task, the
ML approaches utilize existing DTA data, meaning that fewer DTAs overall might
be needed.

A factor affecting the running time of all approaches is the number of demand
variables to be estimated. Increasing the number of variables produces a more
complex estimation problem that results in many more DTAs having to be solved
for an approach to find improving estimates. To cope with this growth, it is
possible to manually adjust the number of variables, e.g., by merging OD pairs
by combining areas into larger zones or by adjusting the number of sub-intervals
of the analysis period in which the demand is estimated. Moreover, arcs equipped
with sensors can be prioritized, and arcs deemed less critical can be filtered
out, making it possible to scale to larger networks.

\section{Conclusions and Future Work} \label{sec:conclusion}

The results presented in Table~\ref{tab:summary_table} and the plots in
Figure~\ref{fig:assignmat_approach}, \ref{fig:spsa_approach},
\ref{fig:knn_approach}, and \ref{fig:fnn_approach} show that the AM-based
approach (method M1) outperforms, in terms of quality of the estimates produced,
the AM-free approach (method M2), as well as the ML approaches (methods M3-M4).
We also note that even if no seed vectors were supplied, the AM-based approach
could find a reasonable estimate for both the ground-truth demands and the
corresponding observed arc counts.

A further advantage of the AM-based approach is that it is easy to apply as no
parameters will have to be set or tuned before using the approach. Moreover,
only two new OF evaluations and thus DTA computations are needed
in each iteration, while in comparison, the AM-free approach needs up to four.
On the other hand, the disadvantage of the AM-based approach is that it might be
harder to include different, more complex types of data into the optimization
problem, as it can be harder to derive analytical expressions for the gradient
of function terms that take into account more complex data types such as road
queuing length, turning ratios, etc.

If we look at the results obtained with the ML approach in
Table~\ref{tab:summary_table} and Figures~\ref{fig:knn_approach} and
\ref{fig:fnn_approach}, we can conclude that the estimated and observed counts
align somewhat well, while the estimated demands and ground-truth demands do not
align at all. This result seems to remain the same independently of which ML
algorithm is used. Based on these results, the ML approach is not good if one
seeks a demand estimate that is sufficiently close to a given seed
vector. However, it can produce more diverse estimates than, e.g., the
AM-based approach and may be useful for other purposes.

Our conclusion is that the gradient-based approaches work the best. However, the
ML methodology does provide a way to utilize data generated through many DTAs,
and it might be a viable option if certain changes are made. In this case, we
suggest possible directions for future research to improve the ML approach:
\begin{itemize}
    % \item Incorporate more problem-specific structural information about the
    % underlying dynamics in the traffic network that ultimately affect the
    % routing of the users in the network (e.g., as it has been done in
    % \cite{oso2019}).

    \item In the optimization step of the ML approach, investigate why only the
    second OF term, relating to the arcs counts, appears to
    have an effect. That is, determine why the first term pertaining to the
    demands is not being optimized and able to guide the search even when
    a seed vector is supplied.

    \item Incorporate additional constraints into the optimization step of the
    ML approach to further restrict the search space and possibly obtain better
    estimates.

    \item Use other types of information, e.g., about the complete paths used by
    users in the network instead of just the arc counts.

    \item Model a different part of the DODE problem and change the type of
    prediction being made (i.e., Eq. \eqref{eqn:objective_function_map}) by
    using more sophisticated ML models. For example, it can be conjectured that
    by using Graph Neural Networks (GNN), it is possible to capture the cause
    and effect better of the demand on the resulting arc counts as the
    relationships and interdependencies that exists between the arcs and thus
    arc counts in the traffic network are modeled more explicitly.
\end{itemize}

% References
\clearpage
\bibliographystyle{tfcad}
\bibliography{bibliography.bib}

% Appendix
\appendix

\section{Derivation of the Gradient of the Objective Function} \label{sec:grad_derivation}

In the following we show how the gradient of the OF $F$ is
derived. To do this we only show how the gradient of the first term is computed,
as the derivation of the gradient of the second term is very similar. The
gradient of the OF $F$ is the defined as:
\begin{align}
    \frac{\partial F}{\partial x_{ws}} = \omega_1 \cdot \frac{\partial f^{(1)}}{\partial x_{ws}} + \omega_2 \cdot \frac{\partial f^{(2)}}{\partial x_{ws}}, \qquad \forall w \in W, s \in S.
\end{align}

To proceed we compute the gradient of the first term $f^{(1)}$ as follows:

\begin{align}
  \frac{\partial f^{(1)}}{\partial x_{ws}} &= \frac{\partial }{\partial x_{ws}} \left(
  \frac{
    \left(
      \sum\limits_{  \substack{w' \in W \\ s' \in S}  } \left(x_{w's'} - \tilde{x}_{w's'} \right)^2
    \right)^{ \nicefrac{1}{2} }
  }
  {
    \left(
      \sum\limits_{  \substack{w' \in W \\ s' \in S}  } \tilde{x}_{w's'}^{2}
    \right)^{ \nicefrac{1}{2} }
  }
  \right) %, \qquad \forall w \in W, s \in S.
\end{align}

We notice that the expression in the denominator is simply a constant. We can factor this constant out and get:

\begin{align}
  = \frac{\partial }{\partial x_{ws}} \left(
      \left(
        \sum\limits_{  \substack{w' \in W \\ s' \in S}  } \left(x_{w's'} - \tilde{x}_{w's'}\right)^2
      \right)^{ \nicefrac{1}{2} }
    \right)
    \cdot 
        \frac{1}{\left( \sum\limits_{  \substack{w' \in W \\ s' \in S}  } \tilde{x}_{w's'}^{2} \right)^{ \nicefrac{1}{2} }}
\end{align}

We then proceed by applying the \emph{chain rule} to this expression and get:

\begin{align}
    =
    \underbracket{
    \frac{
      1
    }
    {
      2
    }
    \cdot \left( \sum\limits_{ \substack{w' \in W \\ s' \in S  }} \left(x_{w's'} - \tilde{x}_{w's'}\right)^2\right)^{-\nicefrac{1}{2}}}_{\textnormal{Deriv. of outer func.}}
    \cdot \underbracket{2 \cdot \left(x_{ws} - \tilde{x}_{ws}\right) 
    \vphantom{\left(\frac{1}{m \cdot n_S} \cdot \sum\limits_{ \substack{w' \in W \\ s' \in S  }} \left(x_{w's'} - \tilde{x}_{w's'}\right)^2\right)^{-\frac{1}{2}}} % Create appropriate vertical space
    }_{\textnormal{Deriv. of inner func.}}
    \cdot \; \frac{1}{\left( \sum\limits_{  \substack{w' \in W \\ s' \in S}  } \tilde{x}_{w's'}^{2}\right)^{\nicefrac{1}{2}}}
\end{align}

When taking the partial derivative of the inner function w.r.t. $x_{ws}$ we notice that the terms with $w' \neq w$ and $s' \neq s$ are $0$. We are thus only left with a single term, when taking the derivative of the inner function. Finally, we can reduce this expression and define the gradient of $f^{(1)}$ as:

\begin{align}
  g^{(1)}_{ws} =
  \frac{ \partial f^{(1)}  }{  \partial x_{ws}  } =
      \frac{
          x_{ws} - \tilde{x}_{ws}
      }
      {
          \left(
              \sum\limits_{  \substack{w' \in W \\ s' \in S  }  } \left(x_{w's'} - \tilde{x}_{w's'}\right)^2
          \right)^{\nicefrac{1}{2}} \cdot \left(\sum\limits_{  \substack{w' \in W \\ s' \in S}  } \tilde{x}_{w's'}^{2}\right)^{\nicefrac{1}{2}}
      }%, \quad \forall w \in W, s \in S
\end{align}

The gradient of the second term $f^{(2)}$, which describes the marginal effect
of a change in the demand variables on the arc counts, can be computed in a
similar way. The only difference is that the relation in
Eq.~\eqref{eqn:assignmat_relation} should be used in the derivation as well.

\section{Machine Learning Algorithms} \label{sec:ml_models}

\textbf{Feed-Forward Networks (FFN):} A feed-forward network is a nonparametric
    model and a type of artificial neural network with a directed acyclic graph
    structure consisting of \emph{units} arranged in layers: First an input
    layer, then possibly several \emph{hidden layers}, and lastly, an
    \emph{output layer}. Each of the units in a layer has links to the units in
    a subsequent layer, and each of the links has an associated weight that
    determines the strength and sign of the connection. In each unit, an
    \emph{activation function} is applied to an offset term (usually referred to
    as a bias term in the FNN terminology) plus the weighted
    sum of activations coming from the units in a previous layer. Depending on
    the layer, different activation functions are used. No activation function
    is used in the input layer and a \emph{ReLU} ($\phi(\cdot) = \max(0,
    \cdot)$) activation function is usually used in the units situated in the
    hidden layers. The choice of activation function used in the output layer is
    tied to the nature of the task at hand. For example, a linear activation
    function is used for the regression task at hand where the target variables
    are continuous. A loss function that measures a model's goodness-of-fit also
    needs to be defined to determine the bias terms and the weights of the links
    between the units. Just as the activation function used in the output layer
    is tied to the task at hand, the choice of the loss function is also closely
    tied to the activation function used in the output layer. For example, a
    Mean Squared Error (MSE) loss function is used for a regression task.

    A FNN architecture is defined by the number of layers (the
    depth of the network $m_d \in \mathbb{N}$ not counting the input layer) and
    the number of units in each layer (the width of a layer $l_i \in
    \mathbb{N}$, with $i \in \left\{1, \dots, m_d \right\}$), along with the
    activation functions used in the different layers. In this case, to train
    the network, the problem is to find a set of weight matrices $\mathcal{W} =
    \left\{ \mathbf{W}_1, \dots, \mathbf{W}_n \right\}$ and corresponding biases
    $\mathcal{B} = \left\{ \mathbf{b}_1, \dots, \mathbf{b}_n \right\}$ that
    minimize the Mean Squared Error (MSE) loss function:
    \begin{align}
        (\hat{\mathcal{W}}, \hat{\mathcal{B}}) &\in \arg\min_{\mathcal{W},\mathcal{B}} ||\mathbf{y} - \hat{\mathbf{y}}_{\mathcal{W}, \mathcal{B}}(\mathbf{x})||^{2}_2 + \lambda \sum^{n}_{i=1}||\mathbf{W}_i||^{2}_2 \label{eqn:nn_opt}\\
        \hat{\mathbf{y}}_{\mathcal{W},\mathcal{B}}(\mathbf{x}) &= \left(h_{m_d} \circ h_{m_d - 1} \hdots \circ h_1 \right)(\mathbf{x}) \label{eqn:nn_estimator}
    \end{align}
    where $h_i(\mathbf{a}) = \phi_{i}\left(\mathbf{W}_{i}^\top \mathbf{a} +
    \mathbf{b}_i \right)$ for $i \in \left\{ 1, \dots, m_d \right\}$.
    Furthermore, the activation function $\phi_i$ is applied \emph{element-wise}
    to the vector-valued argument, and $\mathbf{a}$ contains all the activations
    received in the $i$th layer, while $\mathbf{W}_i$ contains the corresponding
    weights, and $b_i$ is the bias term. 
    Finally, in the last term of Eq.~\eqref{eqn:nn_opt} $\lambda \in
    \mathbb{R}_+$ is a regularization parameter added to reduce overfitting and
    improve the out-of-sample predictive performance. An appropriate value for
    the regularization parameter is usually chosen using out-of-sample
    cross-validation techniques.   
        
    Gradient-based methods are widely used for solving the problem in
    Eq.~\eqref{eqn:nn_opt}. This is because gradients can be computed
    efficiently using the \emph{backpropagation} algorithm (see, e.g.,
    \cite{goo2016} for an in-depth explanation). Essentially, we learn an
    input-output relationship by applying a sequence of semi-affine nonlinear
    transformations, which in terms of the layers in a FNN can
    be understood as applying the composite map in Eq.~\eqref{eqn:nn_estimator}
    to a given input $\mathbf{x}$.

\textbf{$k$-Nearest Neighbors ($k$-NN):} $k$-NN is one of the simplest
    nonparametric models used for regression. To predict the value
    $\hat{\mathbf{y}}$ of a new unseen point $\mathbf{x}$, $k$ needs to be
    specified. Its value specifies the number of points from a training set
    $\mathcal{D}$ closest to $\mathbf{x}$ that should be used to calculate the
    predicted value $\hat{\mathbf{y}}$. Let $N_k(\mathbf{x}) \subseteq
    \mathcal{D}$ be the set of $k$ points closest to $\mathbf{x}$ in terms of
    the Euclidean distance, then the predicted value of $\mathbf{x}$ is
    calculated as the average of the nearest points:
    \begin{align}
        \hat{y}_k(\mathbf{x}) = \frac{1}{k} \cdot \sum_{(\mathbf{x}', \mathbf{y}') \in N_k} \mathbf{y}'.
    \end{align}
    Unlike the FNN and most other ML
    algorithms, no model must be determined, rather the training data are kept
    and continuously used.
\end{document}